\documentclass[10pt]{article}
\usepackage[title,titletoc]{appendix}
\usepackage{lineno}
\usepackage{mathbbold}
\usepackage{}
\usepackage{graphicx}
\usepackage{amsmath,amssymb,amsthm,amsfonts}
\usepackage{amssymb}
\usepackage{color}
\usepackage{amsfonts}
\usepackage{amssymb}
\usepackage{amsthm}
\usepackage{amsmath}
\usepackage{empheq}
\usepackage{indentfirst}
\usepackage{cite}
\usepackage{mathrsfs}
\usepackage{extarrows}
\allowdisplaybreaks[4]
\newtheorem{thm}{Theorem}[section]
\newtheorem{cor}[thm]{Corollary}
\newtheorem{lem}[thm]{Lemma}\newtheorem{prop}[thm]{Proposition}

\theoremstyle{definition}

\theoremstyle{Remark}
\newtheorem{rem}{Remark}[section]

\numberwithin{equation}{section}

\DeclareMathSymbol{\C}{\mathalpha}{AMSb}{"43} \topmargin-.1in

\newcommand{\beq}{\begin{equation}}
\newcommand{\eeq}{\end{equation}}

\newcommand{\lam}{\lambda}
\newcommand{\alp}{\alpha}

\newcommand{\R}{\mathbb{R}}

\newcommand{\m}{\mathrm{m}}

\def\qed{\hfill$\square$\vspace{6pt}}

\newcommand{\bsub}{\begin{subequations}}
\newcommand{\esub}{\end{subequations}$\!$}

\begin{document}
\title{Singular Behavior of an Electrostatic--Elastic Membrane System with an External Pressure}
\author{Yujin Guo \thanks{Wuhan Institute of Physics and Mathematics,
    Chinese Academy of Sciences, Wuhan 430071, P. R. China.
Email: \texttt{yjguo@wipm.ac.cn}. Y. J.  Guo is partially supported by NSFC [grant number 11671394].},
\,   Yanyan Zhang\thanks{Corresponding author. School of Mathematical Sciences, Shanghai Key Laboratory of Pure Mathematics and Mathematical Practice, East China Normal University,
  Shanghai 200241, P.R. China. Email: \texttt{yyzhang@math.ecnu.edu.cn}. Y. Y. Zhang is sponsored by ``Chenguang Program" supported by Shanghai Education Development Foundation and Shanghai Municipal Education Commission [grant number  13CG20]; NSFC [grant number 11431005]; and STCSM  [grant number 18dz2271000].}
    \ and\ Feng Zhou\thanks{School of Mathematical Sciences, Shanghai Key Laboratory of Pure Mathematics and Mathematical Practice, East China Normal University, Shanghai 200241,   P.R. China. Email: \texttt{fzhou@math.ecnu.edu.cn}. F. Zhou is partially supported by NSFC  [grant numbers 11726613, 11431005]; and STCSM [grant number 18dz2271000].}
}

\date{\today}

\smallbreak \maketitle

\begin{abstract}
We analyze nonnegative solutions of the nonlinear elliptic problem $\Delta u=\frac{\lambda f(x)}{u^2}+P$, where $\lambda>0$ and $P\geq0$ are constants, on a bounded domain $\Omega$ of $\mathbb{R}^N$ ($N\geq 1$) with a Dirichlet boundary condition. This equation models an electrostatic--elastic membrane system with an external pressure $P\geq 0$, where $\lambda >0$ denotes the applied voltage. First, we completely address the existence and nonexistence of positive solutions. The classification of all possible singularities at $x=0$ for nonnegative solutions $u(x)$ satisfying $u(0)=0$ is then analyzed for the special case where $\Omega=B_1(0)\subset\mathbb{R}^2$ and $f(x)=|x|^{\alpha}$ with $\alpha > -2$. In particular, we show that for some $\alpha,$  $u(x)$ admits only the ``isotropic" singularity at $x=0$, and otherwise $u(x)$ may admit the ``anisotropic" singularity at $x=0$.  When $u(x)$ admits the  ``isotropic" singularity at $x=0$, the refined singularity of $u(x)$ at $x=0$ is further investigated, depending on whether $P>0$, by applying Fourier analysis.
\end{abstract}
\vskip 0.2truein



\vskip 0.2truein


\section{Introduction}
In this study, we consider nonnegative  solutions of the following singular elliptic equation:

$$\left\{\arraycolsep=1.5pt
\begin {array}{lll}
\Delta u=\displaystyle\frac{\lambda f(x)}{u^2}+P\   \ &\hbox{ in} \,\  \Omega,\\[2mm]
0\leq u\leq1\ &\hbox{ in} \,\  \Omega,\\[2mm]
u=1\ &\hbox{ on} \,\  \partial \Omega,
\end{array}\right.\eqno{(S)_{\lam ,P}}
$$
where $\lam >0$ and $P\geq 0$ are constants,  $0\le f\in L^\infty(\bar \Omega)$ and $f\not\equiv0$ with $ \Omega \subset\R^N$ ($N\ge 1$). If $u>0$ holds in $\Omega$, then $u$ is called a classical solution of $(S)_{\lam ,P}$; However, if $u(x_0)=0$ holds at some $x_0\in \Omega$, then $u$ is called a singular (weak) solution of $(S)_{\lam ,P}$.  When $P=0$, the elliptic problem $(S)_{\lam ,0}$ has been widely investigated over the past few years, see \cite{DG,DY,EGG,GG1,GZ1,PB,YZ} and the references therein. This equation models (cf. \cite{BP,LW}) an electrostatic--elastic membrane system with an external pressure denoted by $P\ge 0$.  This device consists of an elastic membrane suspended over a rigid ground plate, where the normalized distance between the  membrane and the ground plate is described by $u$ in equation $(S)_{\lam ,P}$. When a voltage, represented here by $\lambda$, is applied, the membrane deflects toward the ground plate and a snap-through may occur when
it exceeds a certain critical value $\lambda^*$ (pull-in voltage).
This creates a so-called ``pull-in instability", which greatly
affects the design of many devices. As remarked in \cite{BP}, the study of such systems is therefore important, not only in the field of electro-hydrodynamics, but also in the study of electrostatic actuators and their importance to the design of MEMS devices in which they are used. The permittivity profile $f$ of $(S)_{\lam ,P}$ is assumed to vanish somewhere and satisfies

\begin{equation}\label{permit}
0\leq f\in C^\epsilon(\bar \Omega) \hbox{ for some }\epsilon \in (0,1]
  \hbox{ and } f\not\equiv 0,
\end{equation}
where $ \Omega \subset\R^N$ and $N=1$ or $2$ for an electrostatic--elastic membrane system. We remark that Beckham and Pelesko \cite{BP} recently studied positive solutions of the elliptic problem $(S)_{\lam ,P}$ in certain special domains $\Omega$, where the interesting mathematical structures, including the existence and nonexistence, bifurcation behavior, and stability, of positive solutions were successfully analyzed and computed. For the parabolic version related to $(S)_{\lam ,0},$ the dependence on $f(x)$ of quenching behavior was studied interestingly in \cite{GS}. Moreover, similar equations in which the nonlinear singular term is replaced by $\frac{\lambda f(x)}{u^{\beta}}$ for general cases $\beta>0$ were also studied in the literature, see \cite{G,Ka,LE} and the reference therein.

Stimulated by \cite{BP,CMV,EGG,GG1,PB,YZ}, the main purposes of this study are to address the complete description, in terms of $(\lam, P)$, of the existence and nonexistence of positive solutions $u$ for $(S)_{\lam ,P}$, and the investigate the possible singular behavior at $x=0$ of nonnegative solutions $u$ satisfying $u(0)=0$.
Towards the first purpose, we denote for convenience $0<\Phi\in H^1_0(\Omega)$ to be the unique positive solution of

\begin{equation}\label{2:a}
 -\Delta \Phi =1 \,\ \mbox{in}  \,\ \Omega ;\quad \Phi =0  \,\ \mbox{on}  \,\
\partial\Omega,\quad
\end{equation}
and $P^*:=\frac{1}{\|\Phi\|_\infty}>0.$ Note that $P^*$ depends only on $\Omega.$ By making full use of $\Phi$, the following theorem is concerned with the existence and nonexistence of positive solutions $u$ for $(S)_{\lam ,P}$.

\begin{thm} \label{thm1} Suppose that $f$ satisfies (\ref{permit}).
\begin{enumerate}
\item If $P\geq P^*$, then there is no positive (classical) solution for $(S)_{\lam ,P}$ as soon as $\lam>0$.

\item If $0\leq P<P^*$, then there exists a constant $\lam ^*_P=\lam ^*_P(f,\Omega )$ satisfying
\begin{equation}\label{thm1:a}
 \frac{4}{27(P^*)^2\sup_\Omega f}(P^*-P)^3\le \lam ^*_P\le \frac{|\Omega|-P\int _{\Omega}\Phi dx}{\int _{\Omega}\Phi f dx},\,\ where \,\  P\int _{\Omega}\Phi dx<\frac{P}{P^*}|\Omega|,
\end{equation}
such that
\begin{enumerate}
\item if $0\leq \lam <\lam ^*_P$, there exists at least one positive (classical)  solution for
$(S)_{ \lam ,P}$.
\item if  $\lam >\lam ^*_P$, there is no positive (classical)  solution for $(S)_{\lam ,P}$.
\end{enumerate}
\item The critical constant $\lam ^*_P=\lam ^*_P(f,\Omega )$ is nonincreasing in $P$ for $0\leq P< P^*.$
\end{enumerate}
\end{thm}

We remark that the existence and nonexistence results of Theorem \ref{thm1} are proved for classical solutions by a standard process, see \cite{GG1,PB,YZ}.
The arguments of \cite{EGG,GG1} and the references therein can be actually used to prove in the weak sense that for any fixed $0\leq P<P^*$, there exists a unique minimal (positive) solution $w_{\lam ,P}$ of $(S)_{\lam ,P}$ for any $0<\lam <\lam ^*_P$, which is monotonic strictly in $\lam$. Moreover, $w^*_P(x)=\lim _{\lam \nearrow \lam ^*_P}w_{\lam ,P}(x)$ exists and solves  $(S)_{\lam ^*_P,P}$ uniquely, which is called the extremal solution of $(S)_{\lam ,P}$. Of course, the unique extremal solution $w^*_P$ of $(S)_{\lam ,P}$ may be either regular or singular (in the sense $\|1-w^*_P\|_\infty=1$), which depends on the dimension $N$ and the profile $f$ as well. More precisely, following \cite{EGG,GG1} and the references therein, the following analytic properties of extremal solutions can be further established: if $0<N\le 7$, then the extremal solution $w^*_P$ of $S_{\lam ,P}$ exists and is regular, whereas for $N>7$, such an extremal solution $w^*_P$ of $S_{\lam ,P}$ exists, but $w^*_P$ may be singular, depending on the profile $f$.
Following Theorem \ref{thm1},
the existence and nonexistence of positive solutions for $(S)_{\lam ,P}$ are illustrated by Figure 1.1 below. One can also check that if $\Omega=B_1(0)$ is a ball in $\R^2$, then $\Phi(x)=\frac{1}{4}(1-|x|^2)$ and $P^*=4$, and hence Theorem  \ref{thm1} seems consistent with the numerical observations of \cite{BP}.

 \begin{figure}[!h]
\centering
 {\scalebox{0.37}[0.37]{\includegraphics{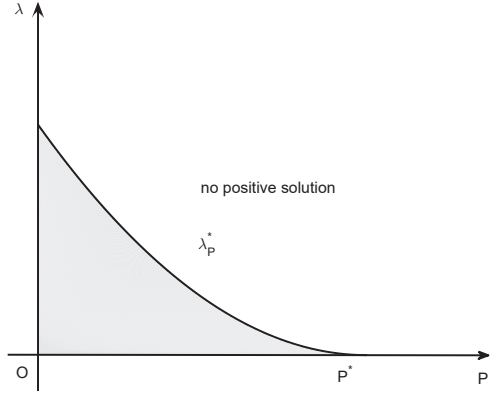}}
 \caption{\label{Fig11}
  {\em  Existence and nonexistence of positive solutions for $(S)_{\lam ,P}$ at different values $(P, \lam)$.}}}
 \end{figure}

The second main purpose of this study is to discuss the singular behavior of
solutions $u$ for $(S)_{\lam ,P}$. This is motivated by the fact that when $f(x)=|x|^4$ and $\Omega=B_1(0) \subset \R^2$, then $u(x)=|x|^2$ is a singular solution of $(S)_{\lam ,P}$ for any $(\lam, P)$ satisfying $\lam +P=4$. In the following, we consider
the equations for $f(x)=|x|^\alp $ ($\alpha >-2$) and $\Omega=B_1(0)\subset\R^2$.
More precisely, we are next concerned with the local behavior near the origin for singular solutions to the following elliptic equation

\begin{equation}\left\{\arraycolsep=1.5pt
\begin {array}{lll}
\Delta u=\displaystyle\frac{\lambda|x|^\alp}{u^2}+P \   \ &\hbox{ in} \,\  B_1(0)\subset\R^2,\\[2mm]
u\geq 0\ \text{and}\ u(0)=0\ &\hbox{ in} \,\  B_1(0),
\end{array}\right.
\label{eqn}
\end{equation}
where $\lam >0$, $\alp>-2$, and $P\geq 0.$

In fact, we rewrite equation \eqref{eqn} as a semilinear evolution elliptic problem

\beq\label{62}
u_{rr}+\frac{1}{r}u_r+\frac{1}{r^2}u_{\theta\theta}=\frac{\lambda r^{\alpha}}{u^2}+P,
\eeq
where $u(x)=u(r, \theta)$, by using the polar coordinate $(r, \theta)\in [0,1]\times S^1$ in $  B_{1}(0) .$
Furthermore, as in \cite{CMV,GZ}, we define $v(t, \theta)$ by

\begin{equation}\label{in:trans}
v(t,\theta ):=r^{-\frac{2+\alp}{3}}u(r,\theta), \ \ \mbox{where}\,\ t=-\ln r \ \ \mbox{and}\ \ r=|x|,
\end{equation}
such that $v(t,\theta)$ satisfies the following evolution elliptic problem

\begin{equation}\label{R:trans-1}
  -v_{tt}+ \displaystyle\frac{2(2+\alp )}{3} v_t=v_{\theta\theta}+ \displaystyle \Big(\frac{2+\alp }{3}\Big)^2 v-\displaystyle\frac{\lam}{v^2}-Pe^{-\frac{4-\alp}{3}t}, \quad (t,\theta)\in (0,+\infty)\times S^1.
\end{equation}
We assume that

\beq\label{78}
\alp\in\mathcal{A} :=\arraycolsep=1.5pt\left\{\begin{array}{lll}
(-2,+\infty),\,\ \mbox{for}\,\ P=0, \\[3mm]
(-2,4),\,\ \mbox{for}\,\ P>0,
\end{array}\right.
\eeq
throughout the remainder of this paper.
In association with the stationary problem of (\ref{R:trans-1}), we also denote $w(\theta)$ to be a solution of

\beq\label{2.6}
w''+\Big({{2+\alpha}\over3}\Big)^2w-\frac{\lambda}{w^2}=0  \ \ \hbox{on}\,\     S^1,
\eeq
and define the solution set by

\beq\label{2.6S}
\mathfrak{S}=\big\{w>0:\, w \ \mbox{is a solution of}\ \eqref{2.6} \big\}.
\eeq
The following analysis of the structure of $\mathfrak{S}$ in terms of $\alp$  plays a fundamental role in classifying the singularities of solutions for \eqref{eqn}.

\begin{thm}\label{thm02}
Consider the set $\mathfrak{S}$ defined in \eqref{2.6S}, where $\lam >0$ is arbitrary. The following results hold:

\begin{enumerate}
\item If

\begin{equation}\label{thm2:3}
\alp \in \mathcal{A}_s :=\arraycolsep=1.5pt\left\{\begin{array}{lll}
\mathcal{A} _0:=
\big(-2,-{1\over2}\big]\cup \displaystyle\bigcup^7_{k=2}\bigg[(k-1)\sqrt 3-2,\frac{3k-4}{2}\bigg],\,\ &\mbox{for}\,\ P=0, \\[4mm]
\mathcal{A} _0\cap(-2,4)=\big(-2,-{1\over2}\big]\cup\big[\sqrt{3}-2,1\big]\cup\big[2\sqrt{3}-2,\frac{5}{2}\big]\cup\big[3\sqrt{3}-2,4\big),\,\ &\mbox{for}\,\ P>0,
\end{array}\right.
\end{equation}
 then we have
 $\mathfrak{S} =\big\{ \big(\frac{9\lam}{(2+\alp )^2}\big)^{\frac{1}{3}}\big\}.$

\item 
If $\alp \in\mathcal{A}\backslash\mathcal{A}_s$ with $\mathcal{A}$ defined in \eqref{78},
     then $\mathfrak{S}$ contains precisely $1+N_0(\alpha)$ connected components:  $\mathfrak{S}_0=\big\{ \Big(\frac{9\lam}{(2+\alp )^2}\Big)^{\frac{1}{3}}\big\}$,  $\mathfrak{S}_1,\cdots ,\mathfrak{S}_i,$ $\cdots ,\mathfrak{S}_{N_0(\alpha)}$. Here,

     \beq \label{2.7}
\mathfrak{S}_i=\big\{w_{j_i}(\cdot+a):\, a\in S^1\big\}, \ i=1,2,\cdots,N_0(\alpha),
\eeq
where $j_i=[\frac{\sqrt{3}(2+\alpha)}{3}]+i,$ $w_{j_i}(\theta)$ satisfying $\min\limits_{\theta \in \R} w_{j_i}(\theta)=w_{j_i}(0)$ is a $\frac{2\pi}{j_i}$-periodic positive solution of \eqref{2.6},  and
$1\leq N_0(\alpha)< +\infty$ denotes the number of integers in $\big(\frac{\sqrt{3}(2+\alpha)}{3},\frac{2(2+\alpha)}{3}\big)$.

\end{enumerate}
\end{thm}

We can get more precise information about case 2 of Theorem \ref{thm02}. In fact, we have
\begin{rem}\label{7474}
Under the assumption of Theorem \ref{thm02},  if $ \alpha \in (1,2\sqrt{3}-2)$, then $\big[\frac{\sqrt{3}(2+\alpha)}{3}\big]=1$, $N_0(\alp)=1,$ $j_1=2,$ $\mathfrak{S}_1=\big\{w_2(\cdot+a):\, a\in S^1\big\}$, and $\mathfrak{S}=\mathfrak{S}_0\cup\mathfrak{S}_1=\big\{ \big(\frac{9\lam}{(2+\alp )^2}\big)^{\frac{1}{3}}\big\}\cup\big\{w_2(\cdot+a):\, a\in S^1\big\}$;
if $ \alpha \in (\frac{5}{2},3\sqrt{3}-2)$, then
$\big[\frac{\sqrt{3}(2+\alpha)}{3}\big]=2$, $N_0(\alp)=1,$ $j_1=3,$ $\mathfrak{S}_1=\big\{w_3(\cdot+a):\, a\in S^1\big\}$, and $\mathfrak{S}=\mathfrak{S}_0\cup\mathfrak{S}_1=\big\{ \big(\frac{9\lam}{(2+\alp )^2}\big)^{\frac{1}{3}}\big\}\cup\big\{w_3(\cdot+a):\, a\in S^1\big\}$, etc.
\end{rem}

We note that for  all $P\ge 0$ and $\alp \geq0$, if $u(x)=u(r,\theta)$ satisfying $u(0)>0$ is a solution of (\ref{eqn}) that is continuous on $B_1(0)\subset \R^2$, then $u\in C^2\big(B_{r^0}(0)\big)$ for some $0<r^0<1.$ This can be obtained by employing the elliptic $L^p$ theory and Schauder theory.
We now return to classify nonnegative solutions in the case where $u(0)=0$, for which we have the following theorem.

%
%
%
%
%
\begin{thm} \label{thm2} Suppose $u(x)=u(r,\theta)$ is a solution of (\ref{eqn}) that is continuous on $B_1(0)\subset \R^2$ and $u(0)=0$. Assume there exist
constants $\beta \in (0,1]$ and $C_{\beta}>0$ such
 that $v$ defined in  (\ref{in:trans}) satisfies the assumption

 \begin{equation}\label{Z:2}
\sup_{0<r\leq1}\frac{1}{r^{\frac{4}{3}+\beta}}\int_{B_{r}(x)}\frac{1}{v}dy \leq C_\beta,\quad \forall \ x=(t,\theta)\in[t^0,+\infty)\times S^1
\end{equation}
for some $t^0\in \R.$ Then, there exists $w\in\mathfrak{S}$ such that

    \beq\label{B65}
\|r^{-\frac{2+\alp}{3}} u(r,\theta)-w(\theta)\|_{C^2(S^1)}\leq C (1-\ln r)^{-\frac{\bbtheta}{7(1-2\bbtheta)}}\quad \text{as}\ \ r\rightarrow 0^+,
\eeq
for some $\bbtheta\in(0,\frac{1}{2})$ depending on $w$.
\end{thm}

We remark that if $\lambda=1,$ $\alpha=1$ and $P=0,$ then $u(x)=|x|$ is a singular solution of \eqref{eqn}, for which case the assumption \eqref{Z:2} holds for $\beta=2/3$ and  $C_{\beta}=\pi.$ By combining Theorem \ref{thm02} with Theorem \ref{thm2} (2), we immediately obtain the following:

\begin{cor}\label{74}
Suppose $u(x)=u(r, \theta)$ is a singular solution of (\ref{eqn}) with $u(0)=0$
and (\ref{Z:2}) holds.
\begin{enumerate}

 \item[$1$.] If 
 $\alpha\in\mathcal{A}_s ,$
 then

\begin{equation}\label{thm2:1}
  \lim _{r\to 0^+}r^{-\frac{2+\alp}{3}} u(r,\theta)=  \Big(\frac{9\lam}{(2+\alp )^2}\Big)^{\frac{1}{3}} \quad \text{in}\ \ C^2(S^1).
\end{equation}

 \item[$2$.] 
     If $\alp \in\mathcal{A}\backslash\mathcal{A}_s,$ then either (\ref{thm2:1}) holds,
 or there exists $w_{j_i}(\theta+a), 1\leq i\leq N_0(\alpha)$ such that

\begin{equation}
 \lim _{r\to 0^+}r^{-\frac{2+\alp}{3}}u(r, \theta)= w_{j_i}(\theta+a)
\end{equation}
in $C^2(S^1)$, where $\mathcal{A}, N_0(\alpha), j_i, w_{j_i}$ are as in Theorem \ref{thm02}.
\end{enumerate}
\end{cor}

According to Remark \ref{7474}, we emphasize that the following more precise information for Corollary \ref{74}(2) can be obtained: if $ \alpha \in (1,2\sqrt{3}-2)$, then either (\ref{thm2:1}) holds or there exists $w_2(\theta+a)$ such that $\lim _{r\to 0^+}r^{-\frac{2+\alp}{3}}u(r, \theta)= w_2(\theta+a);$ if $ \alpha \in (\frac{5}{2},3\sqrt{3}-2)$, then either (\ref{thm2:1}) holds or there exists $w_3(\theta+a)$ such that $\lim _{r\to 0^+}r^{-\frac{2+\alp}{3}}u(r, \theta)= w_3(\theta+a),$ etc.

There are several other remarks concerning the above results.
Firstly, if $\alp \in \mathcal{A}_s$, then Corollary \ref{74} (1) shows that $u$ satisfying $u(0)=0$ admits only the ``{\em isotropic}" singularity at $x=0$ in the sense of (\ref{thm2:1}). However, Corollary \ref{74} (2) implies that $u(x)$ satisfying $u(0)=0$ may admit the ``{\em anisotropic}" singularity at $x=0$ if $\alp \in\mathcal{A}\backslash\mathcal{A}_s$.  Secondly, one can easily check that $u_c=|x|^{\frac{2+\alp}{3}}$, which admits the ``isotropic" singularity at $x=0$, is always a singular solution to \eqref{eqn} with the boundary condition $u|_{\partial B_1(0)}=1$ for $(\lam, P)=\big(\frac{(2+\alpha)^2}{9}, 0\big)$.
Thirdly, we remark that the method of \cite[Lemma 1.6]{CMV} does not work for proving the convergence result in Theorem \ref{thm2}. To overcome this difficulty
for the case $P=0$, where \eqref{R:trans-1} is an autonomous evolution equation, the convergence result of Theorem \ref{thm2} can be followed from \cite{jend} by employing the $\L$ojasiewicz--Simon method. For the case $P>0,$ where \eqref{R:trans-1} is an asymptotically autonomous evolution equation, the method used in \cite{jend} does not work directly. We shall combine the techniques of \cite{JWZ,BenH} with the methods of \cite{jend} to overcome this difficulty for the case $P>0$.

Next, we follow Corollary \ref{74} to further analyze the refined ``isotropic" singularity of $u$ for the case that (\ref{thm2:1}) holds
with
\begin{equation}\label{61}
\alp\in \mathring{\mathcal{A} }:=\mathcal{A} \backslash\{(k-1)\sqrt{3}-2:k=2,3,4,\cdots\},
\end{equation}
where $\mathcal{A} $ is defined in \eqref{78}.
 We set a new transformation

 \begin{equation}\label{1:R:newtrans}
V(t,\theta )=r^{-\frac{2+\alp}{3}}u(r,\theta)-\Big(\frac{9\lam}{(2+\alp )^2}\Big)^{\frac{1}{3}} \quad \mbox{in}\,\ (t_0 ,+\infty)\times S^1, \,\ t=-\ln r \ \ \mbox{and}\ \ r=|x|.
\end{equation}
By carefully analyzing the asymptotic behavior of $V(t,\theta )$ as $t\to +\infty$, the following refined singular behavior is proved in Section 4:

\begin{thm} \label{thm3} Assume $u(x)=u(r, \theta)$
 is a singular solution of (\ref{eqn}) with $u(0)=0$ such that (\ref{thm2:1}) holds. Then we have the following refined singular behavior:
\begin{enumerate}

\item  If $\alp\in \mathring{\mathcal{A} }\backslash\big\{(2\sqrt{3}-2,{3\over2}\sqrt{10}-2]\cup(3\sqrt{3}-2,4)\big\},$ then for both cases $P=0$ and $P>0,$ once $\alp\in \big((k-1)\sqrt 3-2,k\sqrt 3-2), $ there exist $A_k\in \R $ and $\theta _k\in S^1$ such that

\begin{equation}\label{77}
 \lim_{r\to 0} r^{-\frac{\sqrt{9k^2-2(2+\alp )^2}}{3}}\bigg[u(r,\theta)-\Big(\frac{9\lam}{(2+\alp )^2}\Big)^{\frac{1}{3}}r^{\frac{2+\alp}{3}}\bigg]=A_k\sin (k\theta +\theta _k)\quad \mbox{in}\ \, C^2(S^1),
\end{equation}
where $k=1,2,3,4,\cdots$ for $P=0,$ and $k=2,3$ for  $P>0.$
\item   If $\alp\in \big(2(\sqrt 3-1),{3\over2}\sqrt{10}-2\big]\cup\big(3\sqrt 3-2,4\big)$ and $P=0$, then \eqref{77} still holds.

\item If $\alp\in \big(2(\sqrt 3-1),{3\over2}\sqrt{10}-2\big)\cup\big(3\sqrt 3-2,4\big)$ and $P>0$,
then

\begin{equation}
 \lim_{r\to 0} r^{-2}\bigg[u(r,\theta)-\Big(\frac{9\lam}{(2+\alp )^2}\Big)^{\frac{1}{3}}r^{\frac{2+\alp}{3}}\bigg]=\frac{9P}{36+2(2+\alp )^2}\quad \mbox{in}\  \, C^2(S^1).  \label{thm3:1B}
\end{equation}
\item If $\alp={3\over2}\sqrt{10}-2$ and $P>0$, then there exist $A_3\in \R $ and $\theta _3\in S^1$ such that

\begin{equation}
 \lim_{r\to 0} r^{-2}\bigg[u(r,\theta)-\Big(\frac{9\lam}{(2+\alp )^2}\Big)^{\frac{1}{3}}r^{\frac{2+\alp}{3}}\bigg]= A_3 \bigg(\sin (3\theta +\theta _3)+{P\over 9}\bigg)\quad \mbox{in}\  \, C^2(S^1).
\end{equation}
\end{enumerate}

\end{thm}

Under the assumptions of Theorem \ref{thm3}, if $\alp\in \big(2(\sqrt 3-1),{3\over2}\sqrt{10}-2\big)\cup\big(3\sqrt 3-2,4\big)$, then Theorem \ref{thm3} reveals the following refined ``isotropic" singularity on singular solutions $u$ of (\ref{eqn}): for the case $P>0$, $u$ admits the
``{\em strongly isotropic}" singularity at $x=0$, in the sense that for some $\gamma >\frac{2+\alp}{3}$ the limit

\begin{equation}
\lim _{r\to 0}\frac{\big[u(r,\theta)- \big(\frac{9\lam}{(2+\alp )^2}\big)^{\frac{1}{3}}r^{\frac{2+\alp}{3}}\big]}{r^\gamma}\label{thm3:M}
\end{equation}
does not depend on the angle $\theta$; however, for the case $P=0$, $u$   admits the ``{\em weakly isotropic}" singularity at $x=0$, that is, the limit (\ref{thm3:M}) depends on the angle $\theta$ for some $\gamma >\frac{2+\alp}{3}$. Moreover, because $2>\frac{\sqrt{9k^2-2(2+\alp )^2}}{3}$ for 
$\alp\in\big((k-1)\sqrt 3-2,k\sqrt 3-2\big)$ with   $k=3,4,$ Theorem \ref{thm3} implies that the external pressure $P>0$
enhances the convergence speed of singular solutions $u$ tending to 0 as $r\rightarrow 0.$ We finally emphasize that for $f(x)=|x|^\alpha$,
our analysis of Theorems \ref{thm02}--\ref{thm3} holds for the case where $\alpha>-2,$ and while Theorem \ref{thm1} holds only for the case where $\alpha \ge 0.$

The remainder of this paper is organized as follows. In Section 2, we prove Theorem \ref{thm1} on the  existence and nonexistence of positive (classical) solutions for $(S)_{\lam ,P}$. Section 3 is devoted to the proof of Theorems \ref{thm02} and \ref{thm2} on the classification of singular solutions for \eqref{eqn}. In Section 4, we then complete the proof of Theorem \ref{thm3}, which is concerned with the refined singular behavior near the origin of nonnegative solutions $u$ satisfying $u(0)=0$.
 Finally, the proof of Lemma \ref{thm3.1} is given in Appendix A, and a $\L$ojasiewicz--Simon type inequality is established in Appendix B, which plays an important role in the proof of Theorem \ref{thm2}.

\section{Existence and Nonexistence of Positive Solutions}

In this section, we focus on the proof of Theorem \ref{thm1}, which is stimulated by \cite{GG1,PB,YZ} and the references therein for the case $P=0$. We denote $\Phi_P= P \Phi$, where $\Phi$ is the unique positive solution in $H^1_0(\Omega)$ of (\ref{2:a}); $\Phi_P$
satisfies

\begin{equation}\label{2:c}
 0<\Phi_P=P\Phi\le \frac{P}{P^*} \,\ \mbox{in}  \,\ \Omega  ,\,\ \mbox{and}  \,\ \|\Phi_P\|_\infty=\frac{P}{P^*},
\end{equation}
where $P^*:=\frac{1}{\|\Phi\|_\infty}>0$.

 To consider positive solutions of  $(S)_{ \lam ,P},$ by setting $\tilde{u}=1-u$ in $(S)_{ \lam ,P},$ in this section we work on the following equivalent problem for convenience:

$$
\left\{\arraycolsep=1.5pt\begin{array}{lll}\arraycolsep=1.5pt
     -\Delta \tilde{u} &=& \displaystyle\frac{\lambda
f}{(1-\tilde{u})^2}  +P  \quad \mbox{in}\,\, \Omega; \\
   \hfill 0\leq &\tilde{u}&<1  \quad\quad \quad \quad \quad \,   \mbox{in}\,\,  \Omega \,;\\
 \hfill   \tilde{u}   &=& 0  \,\,\,\quad\quad\quad \quad \ \quad  \, \, \mbox{on}\,\,
\partial\Omega.\quad
\end{array}\right.\eqno{(\tilde{S})_{\lam ,P}}$$

\begin{lem} \label{lem:2.1} For any $0\leq P<P^*$, there exists a finite positive constant $\lam ^*_P=\lam ^*_P(f,\Omega )$ satisfying (\ref{thm1:a}) such that

\begin{enumerate}
\item If $0\leq \lam <\lam ^*_P$, there exists at least one solution for
$(\tilde{S})_{ \lam ,P}$.
\item If  $\lam >\lam ^*_P$, there is no solution for $(\tilde{S})_{\lam ,P}$.
\end{enumerate}
\end{lem}

\noindent{\bf Proof.} For any fixed $0\leq P<P^*$, define

\begin{equation}
\lam ^*_P=\lam ^*_P(f,\Omega )=\sup\big\{\lam >0\ |\ (\tilde{S})_{ \lam,P } \ {\rm possesses \ at \ least\
one\ solution} \big\} .\label{2:lamP}
\end{equation}
We first prove that $\lam ^*_P>0$ holds for any $0\leq P<P^*$. It is clear that $\underline{u}\equiv 0$ is a subsolution of $(\tilde{S})_{\lam ,P}$ for any $\lam>0$. To construct a supersolution of $(\tilde{S})_{\lam ,P}$ with  $0\leq P<P^*$, set
$\bar u=s\Phi_{P^*}$ with $\frac{P}{P^*}<s<1,$
such that $\bar u\le s$ in $\Omega$ and $\bar u=0$ on $\partial\Omega$. We then have

\begin{equation}\label{2:d}
-\Delta \bar u=P+\displaystyle\frac{P^*}{\sup _\Omega f}\Big(s-\frac{P}{P^*}\Big)(1-s)^2\frac{\sup _\Omega f}{(1-s)^2}
\ge P+\displaystyle\frac{P^*}{\sup _\Omega f}\Big(s-\frac{P}{P^*}\Big)(1-s)^2\frac{ f(x)}{(1-\bar u)^2} \ \ \mbox{in}  \,\ \Omega.
\end{equation}
Setting

\[
g(s)=\displaystyle\frac{P^*}{\sup _\Omega f}\Big(s-\frac{P}{P^*}\Big)(1-s)^2,\,\ \hbox{where}\,\ \frac{P}{P^*}<s<1,
\]
one can obtain that

\[
g(s)\le g\big(\frac{P^*+2P}{3P^*}\big)= \frac{4}{27(P^*)^2\sup_\Omega f}(P^*-P)^3:=\lam _P.
\]
By taking $s=\frac{P^*+2P}{3P^*},$ this implies from (\ref{2:d}) that

\[
-\Delta \bar u\ge  \frac{\lam _P f}{(1-\bar u)^2}+P \ \ \mbox{in}  \,\ \Omega,
\]
which shows that for any fixed $0\leq P<P^*$, $\bar u=\frac{P^*+2P}{3P^*}\Phi_{P^*}>0$ is a supersolution of $(\tilde{S})_{\lam ,P}$ for $0\le \lam \le \lam _P$. By the method of sub-supersolutions, we conclude that for any fixed $0\leq P<P^*$, there is a solution $\tilde{u}_{\lam ,P}$ of $(\tilde{S})_{
\lam ,P}$ for every $\lam \in (0, \lam _P)$, which implies that $\lam ^*_P\ge \lam _P>0$ holds for any $0<P<P^*$.

We next prove the finiteness of $\lam ^*_P$ for any fixed $0\leq P<P^*$. Suppose that $(\tilde{S})_{
\lam ,P}$  has a solution $\tilde{u}$. Multiplying $(\tilde{S})_{\lam ,P}$
by $\Phi$ and integrating over $\Omega$, we obtain

\[
|\Omega|\ge \displaystyle\int _{\Omega}\tilde{u}dx=-\displaystyle\int _{\Omega}\Phi\Delta \tilde{u}dx=\displaystyle
\int _{\Omega}\frac{\lam  \Phi f }{(1-  \tilde{u})^2}dx+P\displaystyle\int _{\Omega}\Phi dx
\ge \lam  \displaystyle\int _{\Omega}\Phi f dx +P \displaystyle\int _{\Omega}\Phi dx,
\]
which then implies that

\[
\lam ^*_P\le \frac{|\Omega|-P\displaystyle\int _{\Omega}\Phi dx}{\displaystyle\int _{\Omega}\Phi f dx}<+\infty,\,\ \hbox{since} \,\  P\int _{\Omega}\Phi dx<\frac{P}{P^*}|\Omega|<|\Omega|.
\]

For any fixed $0<P<P^*$, because
$\lambda_P^*$ is positive and finite, we choose any $\lam \in (0, \lam ^*_P)$
and use the definition of $\lambda ^*_P$ to find $\bar {\lam}\in
(\lam ,\lam ^*_P)$ such that $(\tilde{S})_{\bar \lambda,P}$ has a solution
$\tilde{u}_{\bar {\lam},p}$ satisfying

$$-\Delta \tilde{u} _{\bar {\lam},P}= \frac{\bar {\lam } f }{(1-\tilde{u}_{\bar
{\lam},P})^2}+P\  \hbox{ in }\ \Omega  ; \:\: 0\leq \tilde{u}_{\bar \lambda, P}<1\ \hbox{ in }\ \Omega
 ; \:\: \tilde{u}_{\bar {\lam}, P}=0 \ \hbox{ on } \ \partial\Omega  ,
$$
which implies that $ -\Delta \tilde{u} _{\bar {\lam},P}\geq   \frac{ \lam   f }{(1-\tilde{u}_{\bar
{\lam},P})^2}+P$ in $ \Omega $. This shows that $\tilde{u}_{\bar {\lam}}$ is a supersolution of $(\tilde{S})_{\lambda,P}$.
Because $\underline{u}\equiv 0$  is a subsolution of $(\tilde{S})_{\lambda,P}$, by the method of sub-supersolutions, we deduce that there is a solution $\tilde{u}_{\lam }$ of $(\tilde{S})_{\lambda,P}$ for every $\lam \in (0, \lam ^*_P)$.
Note from the definition of $\lam ^*_P$
that there is no solution of $(\tilde{S})_{\lambda,P}$ for any $\lam
>\lam ^*_P$. This completes the proof of Lemma \ref{lem:2.1}.
\qed

\begin{lem}\label{lem2.3}
If $P\geq P^*$, then there is no solution for $(\tilde{S})_{\lam ,P}$ as soon as $\lam>0$.
 \end{lem}
\noindent{\bf Proof.}
 On the contrary, suppose there exists $P\geq P^*$ such that $(\tilde{S})_{\lam ,P}$ has a solution $0\le \tilde{u}_{\lam ,P}< 1$ for some $\lam >0$. This implies that

\begin{equation}\label{2.5}
-\Delta (\tilde{u}_{\lam ,P}-\Phi_P)=\frac{\lam f }{(1-\tilde{u}_{\lam ,P})^2}\ge 0 \,\ \mbox{in}  \,\ \Omega ;\quad \tilde{u}_{\lam ,P}-\Phi_P   =0  \,\ \mbox{on}  \,\
\partial\Omega.
\end{equation}
Applying the strong maximum principle to (\ref{2.5}), we then obtain that $\Phi_P< \tilde{u}_{\lam ,P}\le 1 \,\ \mbox{in}  \,\ \Omega,$
which is however a contradiction to the fact that $\|\Phi_P\|_\infty=\frac{P}{P^*}\geq1$ by (\ref{2:c}). This completes the proof of Lemma  \ref{lem2.3}.
\qed

\vskip 0.05truein

\noindent\textbf{Proof of Theorem \ref{thm1}.} We complete the proof of Theorem \ref{thm1} by Lemmas \ref{lem:2.1} and \ref{lem2.3}, and part (3) of Theorem \ref{thm1} can be easily established.
\qed

\section{Classification of Singularities}

This section is concerned with the proof of Theorem \ref{thm02} and Theorem \ref{thm2}, which handle the classification of singular solutions for (\ref{eqn}). We reduce \eqref{eqn} into the semilinear evolution elliptic problem \eqref{R:trans-1}, such that it suffices to analyze the long-time profile of $v(t,\theta)$ for \eqref{R:trans-1} as $t\rightarrow +\infty.$ By using the phase-plane method, we first give the proof of Theorem \ref{thm02} in Subsection \ref{3.1}. The proof of Theorem \ref{thm2} is then completed in Subsection
 \ref{3.3} by employing the theory of infinite dimensional dynamical systems as well as the $\L$ojasiewicz--Simon method.

We start with the following crucial local estimates of singular solutions.

\begin{lem} \label{thm3.1} Suppose that $u$ is a nonnegative singular solution of (\ref{eqn})
satisfying $u(0)=0$ and (\ref{Z:2}), where 
$\alp>-2$ satisfies \eqref{78}. Then there exist constants $0<C_1<C_2<+\infty$ such that

\begin{equation}
 C_1|x|^{\frac{2+\alpha}{3}}\le u(x)\le C_2 |x|^{\frac{2+\alp}{3}}\quad \mbox{as}\ \ |x|\to 0 . \label{4:1}
\end{equation}
\end{lem}
Because Lemma \ref{thm3.1} can be established in a similar way to that in \cite{GZ,HW}, we sketch the proof in Appendix A for simplicity.
Recall that the function $v(t,\theta )$ defined in \eqref{in:trans} satisfies the evolution elliptic equation given by

\beq\label{63}
-v_{tt}+ \displaystyle\frac{2(2+\alp )}{3} v_t=v_{\theta\theta}+ \Big(\frac{2+\alp}{3}\Big)^2 v-
\frac{\lam}{v^2}-Pe^{-\frac{4-\alp}{3}t}, \quad (t,\theta)\in (t_0,+\infty)\times S^1.
\eeq
It follows from Lemma \ref{thm3.1} that $0<C_1\le v\le C_2<+\infty $   uniformly in $(t_0 ,+\infty)\times S^1$ for some $t_0>0$. Therefore,
we need only investigate the long-time behavior of the bounded
solution $v$ of (\ref{63}) as $t\to +\infty$.
The following results give some analytic properties of the evolution equation (\ref{63}).

\begin{lem} \label{lemR.1} Under the assumptions of Theorem \ref{thm2}, the following
results hold:\begin{enumerate}
\item There exists  $\delta\in (0,1)$ such that $v(t,\cdot)$, $v_t(t,\cdot)$, $v_\theta
(t,\cdot)$, $v_{tt}(t,\cdot)$, $v_{t\theta}(t,\cdot)$, $v_{\theta\theta}(t,\cdot)$,
$v_{ttt}(t,\cdot)$, $v_{t\theta\theta}(t,\cdot)$, $v_{tt\theta}(t,\cdot)$, and
$v_{\theta\theta\theta}(t,\cdot)$ all remain bounded in $C^\delta (S^1)$ for any $t\in [t_0,+\infty)$,
where $C^\delta (S^1)$ denotes the usual H\"older continuous space on $S^1$.
\item Both $v_{t}(t,\cdot)$ and $v_{tt}(t,\cdot)$ tend to $0$ in $C^0(S^1)$ as $t\to +\infty$.
\item The ``orbit" $ \mathcal{L}:=\{v(t, \cdot ):\, t\ge t_0\}$ of $v$ is relatively compact in $C^2(S^1)$.
\end{enumerate}
\end{lem}

\noindent{\bf Proof.} We prove this lemma in a similar way to \cite{CMV}. Because all coefficients of (\ref{63}) are bounded and Lemma \ref{thm3.1} gives that $0<C_1\le v\le C_2<+\infty $ holds in $(t_0 ,+\infty)\times S^1$, Lemma \ref{lemR.1}(1) is an immediate consequence of  $L^p$ and the Schauder
estimates for (\ref{63}). Moreover, Lemma \ref{lemR.1}(3) follows by directly applying Lemma \ref{lemR.1}(1) and the Ascoli--Arzel\'a Theorem as well.

As for Lemma \ref{lemR.1}(2), multiplying (\ref{63}) by $v_t$ and integrating it by parts with respect to $\theta$ and $t$, we obtain that

\begin{equation}\label{R:vt}
\int ^{+\infty} _{t_0 }\int_{S^1}v^2_td\theta dt<+\infty ,
\end{equation}
owing to the boundedness of $v$, $v_\theta$, and $v_t$.
Set

\[
k(t)=\int_{S^1}v^2_td\theta .
\]
We thus deduce from the boundedness of $v_t$ and $v_{tt}$ that $k'(t)$ is bounded uniformly on $[t_0,+\infty )$. Therefore, we derive from (\ref{R:vt}) that

\[
k(t)\ge 0,\quad \int ^{+\infty}_{t_0 }k(t)dt<+\infty,
\]
which implies that $k(t)\to 0$ as $t\to +\infty$. It then follows from the boundedness of $v_{t\theta}$ that

\begin{equation}\label{R:vt0}
v_t\to 0\quad \mbox{as}\quad t\to +\infty \quad (\mbox{uniformly in} \ \ \theta \in S^1).
\end{equation}
Therefore, the convergence $v_{tt}\to 0$ as $t\to +\infty$ holds by applying (\ref{R:vt0}) and the
boundedness of $v_{tt}$, $v_{ttt} $, and $v_{tt\theta}$, where the details of the proof are omitted for simplicity.
This completes the proof of Lemma \ref{lemR.1}(2), and we have finished.
\qed

\subsection{Structure of set $\mathfrak{S}$}\label{3.1}

In this subsection, we analyze the set $\mathfrak{S}$ defined in \eqref{2.6S}, where the constant $\frac{2+\alpha}{3}$ is replaced by a generic constant $A$.
Our main results are given by the following proposition, from which Theorem \ref{thm02} can be established immediately.
\begin{prop}\label{thm3.3} Consider the set

\begin{equation} \label{thm3.4:00}
\mathfrak{S}:=\Big\{w(\theta )\in C^2(S^1):\, w''+A^2w-\frac{\lam}{w^2} =0,\quad w>0 \Big\}, \end{equation}
where $A>0,~ \lam  >0$  are given constants. Then we have the following results:

\begin{enumerate}
\item If

\begin{equation} \label{thm3.4:1}
A\in A_c:=\big(0,\frac{1}{2}\big]\cup\bigcup_{k=1}^{6}\big[\frac{k}{\sqrt{3}},\frac{k+1}{2}\big],
\end{equation}
 then
 $\mathfrak{S} =\big\{\big(\frac{\lam}{A^2}\big)^{\frac{1}{3}}\big\}.$

\item If $A\not\in A_c$, then $\mathfrak{S}$ contains precisely $1+N_0(A)$ connected components  $\mathfrak{S}_0=\big\{\big(\frac{\lam}{A^2}\big)^{\frac{1}{3}}\big\}$,  $\mathfrak{S}_1,\cdots ,\mathfrak{S}_i,$ $\cdots ,\mathfrak{S}_{N_0(A)}$, where $\mathfrak{S}_i$ is defined by

    \beq \label{thm3.4:2}
\mathfrak{S}_i=\big\{w_{j_i}(\cdot+a);0\leq a<2\pi\big\},\,\ i=1,2,\cdots,N_0(A).
\eeq
Here, $1\leq N_0(A)< +\infty$ denotes the number of integers in $(\sqrt{3}A,2A)$, $j_i=[\sqrt{3}A]+i,$  and $w_{j_i}(\theta)$ satisfying $\min\limits_{\theta \in \R} w_{j_i}(\theta)=w_{j_i}(0)$ is the $\frac{2\pi}{{j_i}}-$periodic positive solution of

\beq \label{6}
w''+A^2w-\frac{\lambda}{w^2}=0  \ \ \hbox{in}\,\     \mathbb R.
\eeq
\end{enumerate}
\end{prop}

To prove Proposition \ref{thm3.3}, we use the standard phase-plane method (cf. \cite{S2,LZ,Yan,444y}). Note that a first integral of \eqref{6} is given by

\beq \label{11}
(w')^2+A^2w^2+ \frac{2\lambda}{w}=E
\eeq
for some constant $E.$
Define

\[
g(w)=A^2w^2+\frac{2\lambda}{w}, \ \ \hbox{where}\,\   w>0,
\]
such that

\[ \label{15}
g(w)\geq g(w_0):=E_0 = g\bigg(\bigg(\frac{\lambda}{A^2}\bigg)^{1/3}\bigg)=3\lambda^{\frac{2}{3}}A^{\frac{2}{3}}>0,
\]
and

\[
g'(w)<0,\,\ 0<w<w_0;\quad g'(w)>0,\,\ w_0<w<+\infty;
\]
\[
g(w)\rightarrow +\infty \ \ \text{as either} \,\ w\rightarrow 0\,\ \text{or}\,\ w\rightarrow +\infty.
\]
As a consequence, \eqref{6} has nontrivial positive solutions if and only if $E>g(w_0).$ Moreover, it is easy to see that any nontrivial solution of problem \eqref{6} has the following two properties: (i) it is periodic; (ii) if $w(\theta)$ is a solution of \eqref{6}, then $w(\theta+a)$ is also a solution of \eqref{6} for any $a\in\mathbb{R}.$

Suppose now that $w(\theta)$ is a nontrivial positive solution of \eqref{6}. Denote $w_1$ (resp. $w_2$) the minimum (resp. maximum) value of $w(\theta).$ Then $w_1$ and $w_2$ are two roots of

\[
g(w)=E \quad \text{for some}\,\ E>g(w_0),
\]
i.e.,
\beq\label{7}
A^2w_1^2+ \frac{2\lambda}{w_1}=E,\ \
A^2w_2^2+ \frac{2\lambda}{w_2}=E.
\eeq
Therefore, by setting $\tau=\frac{w_2}{w_1}$, we conclude from the above that

\beq\label{9}
w_1^3=\frac{2\lambda}{A^2\tau(1+\tau)}.
\eeq

We can assume without loss of generality that $\theta=0$ is a minimum point of $w(\theta)$ and $\theta =L>0$ is a maximum point of $w(\theta)$, such that $w'(\theta)>0$ for any $\theta\in(0,L).$ Thus, there holds $w'(0)=w'(L)=0,$ where $L>0$ is the half-minimum period of $w.$ Note also from \eqref{11} that

\[
d\theta=\frac{dw}{\sqrt{E-A^2w^2-\lambda\frac{2}{w}}},
\]
which implies that

\beq
L(E)=\int_{w_1}^{w_2}\frac{dw}{\sqrt{E-A^2w^2-\lambda\frac{2}{w}}}.
\eeq
By setting
$y=\frac{w}{w_1}$,  $L(E)$ can be rewritten as

\begin{equation}
L(\tau)           = \int_1^\tau \frac{dy}{\sqrt{\frac{E}{w_1^2}-A^2y^2-\lambda\frac{2}{w_1^3y}}}
           = \frac{1}{A}\int_1^\tau\frac{dy}{\sqrt{1+\tau(1+\tau)-y^2-\frac{1}{y}\tau(1+\tau)}},
\label{14}
\end{equation}
where \eqref{7} and \eqref{9} are used. We next address some analytic properties of $L(\tau)$.

\begin{lem}\label{31}
$L(\tau)$ is continuous on $(1,+\infty)$ and satisfies

\begin{equation}
\label{12}\lim_{\tau\rightarrow 1}L(\tau)=\frac{\pi}{\sqrt{3}A},\quad
\lim_{\tau\rightarrow +\infty}L(\tau)=\frac{\pi}{2A}.
\end{equation}
Moreover, $L(\tau)$ is strictly decreasing in $\tau$.
\end{lem}

\noindent{\bf Proof.} Denote $Q(w)=-A^2w+\frac{\lambda}{w^2}$, and let $w_0=(\frac{\lambda}{A^2})^{\frac{1}{3}}$ be  the unique root of $Q(w).$
By \cite[Lemma 3.2]{S2} we then have

\beq \label{3:49}
L(E)\xlongrightarrow{E\rightarrow E_0 }\frac{\pi}{\sqrt{-Q'(w_0)}}=\frac{\pi}{\sqrt{3}A}.
\eeq
Thus, the first equation of \eqref{12} follows directly from (\ref{3:49}), because $\tau\rightarrow 1$ is equivalent to $E\rightarrow E_0=g(w_0)$. By setting $\xi=\frac{y-1}{\tau-1}$, we rewrite \eqref{14} as

\[
L(\tau) =\frac{1}{A}\int_0^1\frac{\tau-1}{\sqrt{1+\tau(1+\tau)-(\xi(\tau-1)+1)^2-\frac{1}{(\xi(\tau-1)+1)}\tau(1+\tau)}}d\xi.
\]
We then obtain that

\[
\lim_{\tau\rightarrow +\infty}L(\tau)=\frac{1}{A}\int_0^1\frac{d\xi}{\sqrt{1-\xi^2}}=\frac{\pi}{2A},
\]
i.e., the second equation of \eqref{12} holds. Finally, the monotonicity of $L$ in $\tau$ follows directly from the case of $\alpha={1\over2}$ of  \cite[Corollary 5.6]{Ben}.
\qed

\noindent{\bf Proof of Proposition \ref{thm3.3}.}
By Lemma \ref{31}, the range of the half period $L$ is $I_A:=(\frac{\pi}{2A},\frac{\pi}{\sqrt{3}A})={1\over A}(\frac{\pi}{2},\frac{\pi}{\sqrt{3}}).$ Thus \eqref{6} has no nontrivial solution if and only if the interval $I_A$ does not contain $\frac{\pi}{j}$ for any integer $j\geq 1,$ which implies that either $0<A \leq 1/2$ or

\[
\frac{\pi}{j+1}\leq \frac{\pi}{2A}    \quad \text{and}\quad \frac{\pi}{\sqrt{3}A}\leq\frac{\pi}{j}
\quad \text{for some } \ j\geq 1,
\]
i.e., either $0< A \leq 1/2$ or

\beq \label{3:51}
\frac{j}{\sqrt{3}}\leq A\leq \frac{j+1}{2}\quad \text{for some }\ j\geq 1.
\eeq
Note that $\frac{j}{\sqrt{3}}<\frac{j+1}{2}$ holds only for $j\leq 6.$ Thus, \eqref{6} has no nontrivial solution if and only if

$$ A\in A_c:=\big(0,\frac{1}{2}\big]\cup\bigcup_{k=1}^{6}\big[\frac{k}{\sqrt{3}},\frac{k+1}{2}\big].$$
Therefore, if  $ A\in A_c$, we then have $\mathfrak{S} =\big\{\big(\frac{\lam}{A^2}\big)^{\frac{1}{3}}\big\}.$

Suppose now that $A\not\in A_c$. Then it is clear that $\mathfrak{S}_0 =\big\{\big(\frac{\lam}{A^2}\big)^{\frac{1}{3}}\big\}\subset \mathfrak{S}$, and $\mathfrak{S}$ also contains nontrivial solutions because the interval $I_A$ contains $\frac{\pi}{j}$ for some integer $j\geq 1$. In other words, there exists $j\geq1$ such that $\frac{\pi}{2A}<\frac{\pi}{j}<\frac{\pi}{\sqrt{3}A},$ i.e., $\sqrt{3}A<j<2A.$ More precisely,
if we denote by $1\leq N_0(A)<+\infty$ the number of integers in $(\sqrt{3}A,2A),$ then $I_A$ contains
$\{\frac{\pi}{j_i}|j_i=[\sqrt{3}A]+i,i=1,2,\cdots,N_0(A)\}.$
 This implies that \eqref{6} has $N_0(A)$ periodic solutions $w_{j_1}(\theta),w_{j_2}(\theta),\cdots,w_{j_{N_0(A)}}(\theta),$ where each $w_{j_i}(\theta)$ has the $\frac{2\pi}{j_i}-$period and $\min\limits_{\theta} w_{j_i}(\theta)=w_{j_i}(0).$
 Therefore, if $A\not\in A_c$, then $\mathfrak{S}$ contains precisely $1+N_0(A)$ connected components  $\mathfrak{S}_0=\big\{\big(\frac{\lam}{A^2}\big)^{\frac{1}{3}}\big\}$, $\mathfrak{S}_1,\cdots ,\mathfrak{S}_i,$ $\cdots ,\mathfrak{S}_{N_0(A)}$, where $\mathfrak{S}_i$ is defined by (\ref{thm3.4:2}). This completes the proof of the proposition.
\qed

We finally remark that Theorem \ref{thm02} follows immediately from Proposition \ref{thm3.3}, because if

\begin{equation}\label{3.29}\arraycolsep=1.5pt\begin{array}{lll}
 \alpha \in  \displaystyle\big(-2,-{1\over2}\big]\cup \displaystyle\bigcup^7_{k=2}\Big[(k-1)\sqrt 3-2,\frac{3k-4}{2}\Big],
\end{array}\end{equation}
then $0<A:=\frac{2+\alp}{3}\in A_c$, where the set $A_c$ is as in Proposition \ref{thm3.3}.

\subsection{Proof of Theorem \ref{thm2}}\label{3.3}

This subsection is devoted to the proof of Theorem \ref{thm2}, for which we still suppose that $u$ is a singular solution of (\ref{eqn}) satisfying $u(0)=0$ and (\ref{Z:2}). Let $v$ be a solution of (\ref{63}) such that  $0< C_1\le v\le C_2 < \infty$ holds.
We define the ``$\omega$-limit set" $\omega(v)$ of
$v$ by

\begin{equation}\label{R:omega}
\omega(v)=\big\{w\big|\, w\in C^2(S^1),  \;  \exists \;t_n\rightarrow +\infty, \lim_{n\rightarrow+\infty}\|v(t_n,\cdot)-w(\cdot)\|_{C^2(S^1)}=0\big\}.
\end{equation}
A standard argument of dynamical systems then gives that $\omega(v)$ is nonempty, compact, and connected in
$C^2(S^1)$. Note from Lemma \ref{lemR.1} that $\omega(v) \subset \mathfrak{S}$, where $\mathfrak{S}$ is given by \eqref{thm3.4:00}
with $A=\frac{2+\alpha}{3}$ and $\alp >-2$ satisfying \eqref{78}.

Following the above analysis, inspired by \cite{JWZ,jend,BenH}, we further obtain the following convergence result.
\begin{prop}\label{B64}
Under the assumption of Theorem \ref{thm2}, let $0<C_1\le v\le C_2<\infty$ be a solution of the evolution equation

\beq\label{B39}
-v_{tt}+ 2A v_t=v_{\theta\theta}+ A^2v-\displaystyle\frac{\lam}{v^2}-Pe^{-\beta t}, \quad (t,\theta)\in (t_0,+\infty)\times S^1,
\eeq
where $A, \beta,\lam>0, P\geq 0$ are given constants. Then there exists a positive solution $w$ of

\begin{equation}
 w_{\theta\theta}+ A^2w-\frac{\lam}{w^2} =0 \ \ in\,\ S^1,
 \end{equation}
such that

 \beq \label{B63}
\|v(t,\cdot)-w(\cdot)\|_{C^2(S^1)}\leq C(1+t)^{-\frac{\bbtheta}{7(1-2\bbtheta)}}\,\ \mbox{as}\,\ t\to\infty,
 \eeq
where $\bbtheta\in(0,\frac{1}{2})$ is a constant depending on $w$.
\end{prop}

In order to prove Proposition \ref{B64}, we need to borrow the following technical lemma, which was established in \cite{FS,HT}:

\begin{lem}\label{B55}
 Let $0\le \mathcal {Z}\in L^2\big((t_0,+\infty)\big)$ be a measurable function on $(t_0,+\infty)$ and $\zeta\in (0,\frac{1}{2})$. If there exist two constants $C>0$ and $T_0\geq t_0$ such that

\[
\int_t^{+\infty}\mathcal {Z}^2(s)ds\leq C \mathcal {Z}^{\frac{1}{1-\zeta}}(t) \quad \text{for a.e. }\ t\geq T_0,
\]
then $\mathcal {Z}\in L^1(T_0,+\infty).$
\end{lem}

\noindent\textbf{Proof of Proposition \ref{B64}.} Because $\omega(v)$ defined by \eqref{R:omega} is a nonempty, compact, and connected subset of $\mathfrak{S}$, we take
$w\in\omega(v)$ and a sequence $\{t_n\}$ such that

\beq \label{518}
v(t_n,\cdot)\rightarrow w\quad\text{as}\,\ t_n\rightarrow +\infty
\eeq
in $C^2(S^1)$. For convenience, we denote

\beq\label{70}
j(v)=A^2v-{\lambda\over v^2},\quad z(t)=Pe^{-\beta t}.
\eeq
In the following, we shall prove that $\omega(v)$ contains a single element $w$, and it satisfies the estimate \eqref{B63} for large $t$.
The proof is divided into the following four steps:

\vskip 0.05truein

{\bf  Step 1:} For any $\varepsilon>0$, define for all $t\geq t_0$,

\beq\label{B41}
H(v)=-{1\over2}\int_{S^1}|v_t|^2d\theta+(1+2A\varepsilon) E(v)+\varepsilon(v_{\theta\theta}+j(v),v_t),
\eeq
where

\[
E(v)=\int_{S^1}\Big(\frac{1}{2}v_{\theta}^2-J(v)\Big)d\theta,\,\ J(v)=\frac{A^2}{2}v^2+\frac{\lambda}{v}.
\]
We claim that for $\varepsilon>0$ small enough,

\beq\label{67}
H(v)\equiv H_{\infty}\quad \text{and}\quad E(v)\equiv E_{\infty}\quad \text{on} \ \omega(v),
\eeq
 where $H_{\infty}$ and $E_{\infty}$ are two constants depending on $\varepsilon$.

To prove the above claim, we first note that \eqref{B39} can be rewritten as

\beq\label{B40}
-v_{tt}+ 2A v_t=v_{\theta\theta}+ j(v)-z(t), \quad (t,\theta)\in (t_0,+\infty)\times S^1,
\eeq
by
(\ref{70}).
Multiplying \eqref{B40} by $v_t$ and integrating on $S^1$, we obtain that

\beq\label{B42}
\frac{d}{dt}\bigg(-\frac{1}{2}\|v_t\|_{L^2(S^1)}^2+E(v)\bigg)=-2A\|v_t\|_{L^2(S^1)}^2-\int_{S^1}z(t)v_td\theta,
\eeq
which implies that

\beq\label{B43}
\frac{dH}{dt}=-2A\|v_t\|_{L^2(S^1)}^2-\int_{S^1}z(t)v_td\theta+2A\varepsilon \frac{dE}{dt}+\varepsilon(v_{\theta\theta}+j(v),v_t)_t.
\eeq
By \eqref{B40}, we also have

\begin{equation}\label{B45}
\begin{split}
(v_{\theta\theta}+j(v),v_t)_t
=&-\|v_{\theta t}\|_{L^2(S^1)}^2+\int_{S^1}j'(v)v_t^2d\theta -\|v_{\theta\theta}+j(v)\|_{L^2(S^1)}^2+(v_{\theta\theta}+j(v),z(t))\\
&+4A^2\|v_t\|_{L^2(S^1)}^2
+2A(z(t),v_t)
-(v_{tt},2Av_t).
\end{split}
\end{equation}
Together with \eqref{B43}, this yields

\begin{equation}\label{B48}
\begin{split}
\frac{dH}{dt}
=&\int_{S^1}(-2A+\varepsilon j'(v))v_t^2d\theta-\varepsilon\|v_{\theta\theta}+j(v)\|_{L^2(S^1)}^2-\varepsilon\|v_{\theta t}\|_{L^2(S^1)}^2\\
&-(z(t),v_t)+\varepsilon(v_{\theta\theta}+j(v),z(t))\\
\leq&\int_{S^1}(-2A+\varepsilon j'(v)+\varepsilon)v_t^2d\theta-\frac{\varepsilon}{2}\|v_{\theta\theta}+j(v)\|_{L^2(S^1)}^2-\varepsilon\|v_{\theta t}\|_{L^2(S^1)}^2+C\|z(t)\|_{L^2(S^1)}^2.
\end{split}
\end{equation}
Therefore, there exists a constant $K>0$ and $\varepsilon_0 >0$ such that for any $\varepsilon \in (0,\varepsilon_0)$,
the following holds

\beq\label{B46}
{dH\over{dt}}-K\|z(t)\|_{L^2(S^1)}^2\leq 0.
\eeq
That is, $\frac{d}{dt}\tilde{H}(t)\leq 0$, where $\tilde{H}(t)=H(v(t))+K\int_t^{+\infty}\|z(s)\|_{L^2(S^1)}^2ds.$
Because $\tilde{H}$ is bounded from below by Lemmas \ref{thm3.1} and \ref{lemR.1}, we infer that for any $\varepsilon \in (0,\varepsilon_0)$, $\tilde{H}\rightarrow H_{\infty}$ as $t\rightarrow +\infty$ for some constant $H_{\infty}$. Because

\[
\lim_{t\rightarrow +\infty}\int_t^{+\infty}\|z(s)\|_{L^2(S^1)}^2ds=\lim_{t\rightarrow +\infty}2\pi P^2\int_t^{+\infty}e^{-2\beta s}ds=0,
\]
we obtain that

\beq\label{69}
\lim_{t\rightarrow +\infty}H(v(t))=H_{\infty},
\eeq
and hence

\[
\lim_{t\rightarrow +\infty}E(v(t))=\frac{1}{1+2A\varepsilon}\lim_{t\rightarrow +\infty}H(v(t))=\frac{1}{1+2A\varepsilon}H_{\infty}:=E_{\infty},
\]
where Lemma \ref{lemR.1} is used.
By the definition of $\omega(v),$ it is easy to see that $E(v)\equiv E_{\infty}$, which implies that
$H(v)\equiv H_{\infty}$ on $\omega(v)$, and the claim is therefore proved.

\vskip 0.05truein

{\bf  Step 2:} We claim that there exist $\bbtheta \in(0,\frac{1}{2})$ and $T_1>0$ such that

\beq \label{520}
|H(v)-H_{\infty}|^{1-\bbtheta}
\leq  C\big(\|v_t\|_{L^2(S^1)}+\|v_{\theta\theta}+j(v)\|_{L^2(S^1)}\big),\quad \text{for all}\ t>T_1.
\eeq
To prove (\ref{520}), we first note from Lemma \ref{516} and Step 1 that for each $v_{\infty}\in \omega(v),$ there exist constants $\sigma_{v_{\infty}}>0$ and $\bbtheta_{v_{\infty}}\in(0,\frac{1}{2})$ depending on $v_{\infty}$, such that

\beq\label{T8}
\big\|v_{\theta\theta}+j(v)\big\|_{L^2(S^1)}\geq |E(v)-E(v_{\infty})|^{1-\bbtheta_{v_{\infty}}}=|E(v)-E_{\infty}|^{1-\bbtheta_{v_{\infty}}},\,\ v\in\mathcal {B}_{\sigma_{v_{\infty}}}(v_{\infty}),
\eeq
where we denote the set

\[
\mathcal {B}_{\sigma_{v_{\infty}}}(v_{\infty}):=\big\{v\in C^2(S^1):\|v-v_{\infty}\|_{C^2(S^1)}<\sigma_{v_{\infty}}\big\}.
\]
Because the union of balls $\{\mathcal {B}_{\sigma_{v_{\infty}}}(v_{\infty}):v_{\infty}\in \omega(v)\}$ forms an open cover of $\omega(v)$, by the compactness of $\omega(v)$ in $C^2(S^1),$ there exist $v^i_{\infty}\in \omega(v)$ $(i=1, 2, \cdots, m)$ such that $\cup_{i=1}^{i=m}\mathcal {B}_{\sigma_i}(v^i_{\infty})$ $(i=1, 2, \cdots, m)$ is a subcover of
 $\omega(v),$
 where the constants $\sigma_i:=\sigma_{v^i_{\infty}}$ and $\bbtheta_i=\bbtheta_{v^i_{\infty}}$ corresponding to $v^i_{\infty}$ are as in (\ref{T8}).
From the definition of $\omega(v),$ there exists a sufficiently large $T_0>t_0$ such that

\[
v(t)\in \bigcup_{i=1}^{m}\mathcal {B}_{\sigma_i}(v^i_{\infty}),\, \ t\geq T_0.
\]
Because $v^i_{\infty}\in \omega(v)\subset \mathfrak{S},$ by taking

\beq\label{68}
\bbtheta=\min\{\bbtheta_i, i=1, 2, \cdots, m\}\in\big(0,\frac{1}{2}\big),
\eeq
we deduce from \eqref{67} and (\ref{T8}) that

\beq \label{B47}
\big\|v_{\theta\theta}+j(v)\big\|_{L^2(S^1)}\geq |E(v(t))-E_{\infty}|^{1-\bbtheta},  \quad   t\geq T_0.
\eeq

Using the H\"older inequality, we obtain from (\ref{B41}) that for any $w\in\omega(v),$

\beq \label{517}
|H(v)-H(w)|^{1-\bbtheta}\leq C_1\big(\|v_t\|_{L^2(S^1)}^{2(1-\bbtheta)}+|E(v)-E(w)|^{1-\bbtheta}+
\|v_{\theta\theta}+j(v)\|_{L^2(S^1)}^{1-\bbtheta}\|v_t\|_{L^2(S^1)}^{1-\bbtheta}\big)
\eeq
holds for some constant $C_1>0$.
Because Young's inequality yields that

\[
\|v_{\theta\theta}+j(v)\|_{L^2(S^1)}^{1-\bbtheta}\|v_t\|_{L^2(S^1)}^{1-\bbtheta}\leq\|v_{\theta\theta}+j(v)\|_{L^2(S^1)}
+C\|v_t\|_{L^2(S^1)}^{1-\bbtheta\over \bbtheta},
\]
we obtain from \eqref{517} that

\beq \label{517H}
|H(v)-H(w)|^{1-\bbtheta}\leq C\big(\|v_t\|_{L^2(S^1)}^{2(1-\bbtheta)}+|E(v)-E(w)|^{1-\bbtheta}+
\|v_{\theta\theta}+j(v)\|_{L^2(S^1)}+\|v_t\|_{L^2(S^1)}^{1-\bbtheta\over \bbtheta}\big).
\eeq
Recall from Lemma \ref{lemR.1} that $\|v_t\|_{L^2(S^1)}\rightarrow 0$ as $t\to\infty$. Because ${1-\bbtheta\over \bbtheta}>1$ and $2(1-\bbtheta)>1$, we conclude from \eqref{67}, \eqref{B47}, and (\ref{517H}) that there exist $T_1>T_0$ and $C>0$ such that (\ref{520}) holds for all $t>T_1$, and Step 2 is therefore proved.

\vskip 0.05truein

{\bf  Step 3:}  We claim that

\beq \label{517MM}
 \big\|v(t)- w\big\|_{C^2(S^1)}\to 0  \ \ \text{as} \ \ t\to\infty,
\eeq
 which implies that $\omega(v)$ contains a single element $w,$ where $w$ is as in \eqref{518}.

Denote

 \beq\label{B61}
\mathcal {Y}(t)=\|v_t\|_{L^2(S^1)}+\|v_{\theta\theta}+j(v)\|_{L^2(S^1)}.
\eeq
Note first from \eqref{B48} that

\beq \label{B49}
\frac{dH}{dt}+C_1\mathcal {Y}^2(t)\leq C_2 \|z(t)\|_{L^2(S^1)}^2.
\eeq
Integrating \eqref{B49} over $(t,\infty)$, where $t>T_1$, we obtain from \eqref{70} and \eqref{69} that

\[
H_{\infty}-H(t)+C_1 \int_t^{+\infty}\mathcal {Y}^2(s)ds\leq C_2\int_t^{+\infty}\|z(s)\|_{L^2(S^1)}^2ds=C_3e^{-2\beta t},
\]
that is,

\beq \label{B50}
C_1 \int_t^{+\infty}\mathcal {Y}^2(s)ds \leq H(t)-H_{\infty} + C_3e^{-2\beta t}.
\eeq
Because it follows from \eqref{520} that

\beq \label{B5}
H(t)-H_{\infty}\leq C\mathcal {Y}^{\frac{1}{1-\bbtheta}}(t),\quad  t>T_1,
\eeq
we then have

\beq \label{B51}
C_1 \int_t^{+\infty}\mathcal {Y}^2(s)ds\leq C\mathcal {Y}^{\frac{1}{1-\bbtheta}}(t)+C_3e^{-2\beta t},\quad  t>T_1.
\eeq
By noting $0<\bbtheta<\frac{1}{2}$, there exists $T_2>T_1$ such that

\beq \label{B52}
\int_t^{+\infty}e^{-4\beta(1-\bbtheta)s}ds=C e^{-4\beta(1-\bbtheta)t}\leq Ce^{-2\beta t},\quad  t>T_2.
\eeq
Define

$$\mathcal {Z}(t)=\mathcal {Y}(t)+e^{-2\beta(1-\bbtheta)t}.$$
We then deduce from \eqref{B51} and \eqref{B52} that

\beq \label{B53}
\int_t^{+\infty}\mathcal {Z}^2(s)ds\leq C\Big(\int_t^{+\infty}\mathcal {Y}^2(s)ds+\int_t^{+\infty}e^{-4\beta(1-\bbtheta)s}ds\Big)
\leq  C \mathcal {Z}^{\frac{1}{1-\bbtheta}}(t),\,\ t>T_2.
\eeq
Applying Lemma \ref{B55}, we thus conclude from \eqref{B53} that

\[
\int_{T_2}^{+\infty}\mathcal {Z}(t)dt< +\infty,
\]
which further implies that

\beq \label{B54}
\int_{T_2}^{+\infty}\|v_t\|_{L^2(S^1)}dt< +\infty.
\eeq
Because

\[
\|v(t)-v(s)\|_{L^2(S^1)}\leq\int_s^t\|v_t(\tau)\|_{L^2(S^1)}d\tau,
\]
we obtain from \eqref{B54}  that

 \beq\label{529}
 v(t)\rightarrow w\quad \text{in} \,\ L^2(S^1)\,\ \text{as} \,\ t\to\infty,
 \eeq
where $w$ is the same as that of \eqref{518}. By the relative compactness of the orbit $\{v(t, \cdot ):\, t\ge t_0\},$ we obtain the desired conclusion \eqref{517MM}.

\vskip 0.05truein

{\bf  Step 4:} We proceed to prove that the convergence rate of \eqref{B63}  holds true.
Essentially, combining \eqref{520} with \eqref{B49}  yields that

\beq\label{66}
\frac{d}{dt}[H(v)-H_{\infty}]+C_1[H(v)-H_{\infty})]^{2(1-\bbtheta)}\leq C_2 e^{-2\beta t},\,\ t>T_1,
\eeq
where $0<\bbtheta<\frac{1}{2}$ is as in Step 2.
By \eqref{66}, direct  calculations give that

\beq\label{B62}
 H(v)-H_{\infty} \leq C (1+t)^{-\frac{1}{1-2\bbtheta}} \,\ \text{for sufficiently large}\ t>0.
\eeq
We then infer from \eqref{B48} and \eqref{B62} that for sufficiently large $t>0$,

\[
\int_t^{2t}\mathcal {Y}^2(s)ds\leq \tilde{H}(t)-\tilde{H}(2t)\leq \tilde{H}(t)-H_{\infty}\leq H(t)+C e^{-2\beta t}-H_{\infty}\leq C (1+t)^{-\frac{1}{1-2\bbtheta}},
\]
which yields that

\[
\int_t^{2t}\mathcal {Y}(s)ds\leq t^{\frac{1}{2}}\Big(\int_t^{2t}\mathcal {Y}^2(s)ds\Big)^{\frac{1}{2}}\leq C(1+t)^{-\frac{\bbtheta}{1-2\bbtheta}}.
\]
We thus have

\[
\int_t^{+\infty}\mathcal {Y}(s)ds\leq \sum_{j=0}^{+\infty}\int_{2^jt}^{2^{j+1}t}\mathcal {Y}(s)ds\leq C\sum_{j=0}^{+\infty}(2^jt)^{-\frac{\bbtheta}{1-2\bbtheta}}\leq   C(1+t)^{-\frac{\bbtheta}{1-2\bbtheta}},
\]
which implies that

\[
 \|v(t)-w\|_{L^2(S^1)}\leq\int_t^{+\infty}\|v_t\|_{L^2(S^1)}ds\leq \int_t^{+\infty}\mathcal {Y}(s)ds\leq   C(1+t)^{-\frac{\bbtheta}{1-2\bbtheta}}.
\]
Using the Sobolev imbedding theorem and Gagliardo--Nirenberg inequality, by Lemma \ref{lemR.1}, we thus conclude that for sufficiently large $t>0$,

\[
\begin{split}
\|v(t)-w\|_{C^2(S^1)} \leq &C \|v(t)-w\|_{H^3(S^1)}\\
\leq &C_1 \|D^3(v(t)-w)\|^{\frac{6}{7}}_{L^{\infty}(S^1)}\|v(t)-w\|_{L^2(S^1)}^{\frac{1}{7}}+C_2\|v(t)-w\|_{L^2(S^1)}\\
\leq & C(1+t)^{-\frac{\bbtheta}{7(1-2\bbtheta)}},
\end{split}
\]
and the estimate \eqref{B63} is then proved. The proof of this proposition is therefore completed.
\qed

By applying Proposition \ref{B64}, we are ready to complete the proof of Theorem \ref{thm2}.

\vskip 0.05truein



\noindent{\bf Proof of Theorem \ref{thm2}.} Under the assumption (\ref{Z:2}), Lemma \ref{thm3.1} is then applicable. In \eqref{B39}, denote $A:=\frac{2+\alp}{3}>0$ for the case where $P=0,$ and further choose $\beta:=\frac{4-\alpha}{3}$ for the  case where $P>0$ and $-2<\alp <4$. Then the convergence of Theorem \ref{thm2}(2) follows directly  from Proposition \ref{B64} in view of \eqref{in:trans}. This completes the proof of Theorem \ref{thm2}.
 \qed

\section{Refined Singular Behavior}

In this section, we prove Theorem \ref{thm3} on the refined singular behavior of solutions $u$ satisfying $\lim _{r\to 0^+}r^{-\frac{2+\alp}{3}} u(r,\theta)=\Big(\frac{9\lam}{(2+\alp )^2}\Big)^{\frac{1}{3}}.$
Throughout this entire section, we define

\begin{equation}\label{R:mu}
\m := \Big(\frac{9\lam}{(2+\alp )^2}\Big)^{\frac{1}{3}}>0,\,\ \mu :=\frac{2+\alp}{3}>0,
\end{equation}
and
always suppose that $\alp$ satisfies \eqref{61}.
 We use the transformation

\begin{equation}\label{R:newtrans}
V(t,\theta )=r^{-\frac{2+\alp}{3}}u(r,\theta)-\m.
\end{equation}
Therefore, $\lim _{t\to +\infty}V(t,\cdot )=0$ and $V(t,\theta )$ is a uniformly bounded solution of the following evolution elliptic equation

\begin{equation}\label{R:3.17'}
 -V_{tt}+ 2\mu V_t=V_{\theta\theta}+ \mu ^2V+\displaystyle\frac{\lam V(V+2\m)}{\m ^2(V+\m)^2}-Pe^{-(2-\mu )t}, \quad (t,\theta)\in (t_0,+\infty)\times S^1,
\end{equation}
where $\mu>0$ and $\m>0$ are as in (\ref{R:mu}). In the following, we investigate the asymptotic behavior of the Fourier coefficients of $V(t,\theta )$ satisfying (\ref{R:3.17'}). We start with the following exponential decay of $V(t,\cdot )$ as $t\to+\infty$.

\begin{lem} \label{R:C} Under the assumptions of Theorem \ref{thm3}, suppose that
$\alp $ satisfies (\ref{61}). Then there exists some constant $\varepsilon >0$ such that

\begin{equation}\label{R:2.1}
\sup _{t\ge t_0}e^{\varepsilon t}\|V(t,\cdot )\|_{C^0(S^1)}<+\infty.
\end{equation}
\end{lem}

\noindent{\bf Proof.} Inspired by \cite{CMV}, on the contrary, suppose that (\ref{R:2.1}) is false. Set $\rho (t)=\|V(t,\cdot )\|_{C^0(S^1)}$; then $\rho (t)\in C^0\big([t_0, +\infty)\big)$, and

\begin{equation}\label{R:2.2} \lim _{t\to +\infty}\rho (t)=0,\quad \lim _{t\to +\infty}
\sup e^{\varepsilon t}\rho (t)=+\infty
\end{equation}
for any constant $\varepsilon >0$. By applying \cite[Lemma A.1]{CMV}, there exists a function $\eta (t)\in
C^\infty\big([t_0, +\infty)\big)$ such that

\bsub \label{4in:1}
\begin{align}
\eta (t)>0,\quad \eta '(t)<0,\quad \lim _{t\to +\infty}\eta (t)=0,\quad \lim _{t\to +\infty} e^{\varepsilon t} \eta (t)=+\infty, \ \text{for any}\ \varepsilon>0, \label{R:2.4}\\
 0<\lim _{t\to +\infty}\sup \frac{\rho (t)}{\eta (t)}<+\infty ,\qquad \qquad \label{R:2.5}\\
 \big(\eta '/\eta \big)',\ \big(\eta ''/\eta \big)'
\in L^1\big((t_0, +\infty)\big),\quad
 \lim _{t\to +\infty} \frac{\eta '(t)}{\eta (t)}=
\lim _{t\to +\infty} \frac{\eta ''(t)}{\eta (t)}=0 .\label{R:2.8}
\end{align}\esub
Define $w(t,\theta )=\frac{V(t,\theta )}{\eta (t)}$, such that $w$ is  bounded uniformly in $[t_0,+\infty)\times S^1$.
Then by (\ref{R:3.17'}), we have $w$ satisfies

\begin{equation}\label{R:2.9}
-w_{tt}+2\Big(\mu -\frac{\eta ' }{\eta  }\Big)w_t=w_{\theta\theta}-P\frac{e^{-(2-\mu )t}}{\eta}
+\Big[\mu ^2+\frac{\eta ''-2\mu \eta '}{\eta }+\frac{\lam (2\m +w\eta )}{\m ^2(\m +w\eta )^2}\Big]w,
\end{equation}
where $(t,\theta)\in (t_0 ,+\infty)\times S^1$. Note from (\ref{R:2.4})--(\ref{R:2.8}) that all coefficients of equation (\ref{R:2.9}) are bounded  uniformly in $(t_0 ,+\infty)\times S^1$, and $\lim _{t\to+\infty}P\frac{e^{-(2-\mu )t}}{\eta(t)}=0$ in view of assumption (\ref{61}). By applying $L^p$ and the Schauder estimates of (\ref{R:2.9}), similar to Lemma \ref{lemR.1}(1) one can deduce that there exists $\delta\in
(0,1) $ such that $w(t,\cdot)$, $w_t(t,\cdot)$, $w_\theta (t,\cdot)$, $w_{tt}(t,\cdot)$,
$w_{t\theta}(t,\cdot)$, $w_{\theta\theta}(t,\cdot)$, $w_{ttt}(t,\cdot)$,
$w_{t\theta\theta}(t,\cdot)$, $w_{tt\theta}(t,\cdot)$, and $w_{\theta\theta\theta}(t,\cdot)$ all
remain bounded in $C^\delta (S^1)$ for all $t\in [t_0,+\infty)$. Applying (\ref{R:2.8}), as in   Lemma \ref{lemR.1}(2), one can further prove that
$w_{t}(t,\cdot)$ and $w_{tt}(t,\cdot)$ tend to $0$ in $C^0(S^1)$-topology as $t\to +\infty$. So if we
define the ``$\omega$-limit set" $\mathrm{\Gamma}(\mathcal{L}')$ of the ``orbit" $\mathcal{L}':=\{w(t,
\cdot ):\, t\ge t_0\}$ for (\ref{R:2.9}) as

\[
\mathrm{\Gamma}(\mathcal{L}'):=  \bigcap _{t\ge t_0}\overline{\Big\{w(\iota ,\cdot );\, \iota \ge t\Big\}} ,
\]
where the closure is with respect to the topology of $C^2(S^1)$, then a standard argument of dynamical systems shows that $\mathrm{\Gamma}(\mathcal{L}')$ is a nonempty, compact, and
connected set in $C^2(S^1)$. Moreover, $\mathrm{\Gamma}(\mathcal{L}') \subset \mathfrak{S}'$, where
$\mathfrak{S}'$ is the set of stationary solutions of (\ref{R:2.9}), i.e.,

\[
\mathfrak{S}':=\Big\{w(\theta )\in C^2(S^1):\, w''+3\mu ^2w=0 \Big\}.
\]
Because $\sqrt{3}\mu>0$ is not an integer for $\alp$ satisfying \eqref{61}, we obtain $\mathfrak{S}'=\{0\},$ which contradicts \eqref{R:2.5}. This completes the proof of this lemma.
 \qed

The following Fourier analysis gives better estimates of the power $\varepsilon$ in (\ref{R:2.1}), depending on the specific range of $\alp$ and $P$.

\begin{lem} \label{R:D} Under the assumptions of Theorem \ref{thm3}, suppose $V(t,\theta)$ satisfies (\ref{R:newtrans}) and $\mu$ is given by (\ref{R:mu}).
Then there exists a constant $M>0$ such that
\begin{enumerate}
\item If $\alp\in \mathring{\mathcal{A} }\backslash\big\{(2\sqrt{3}-2,{3\over2}\sqrt{10}-2]\cup(3\sqrt{3}-2,4)\big\},$ then for both cases $P=0$ and $P>0,$ there holds
\begin{equation}\label{R:2.1p4}
 \|V(t,\cdot )\|_{C^0(S^1)}\le Me^{-\big(\sqrt{k^2-2\mu ^2}-\mu \big) t},\quad \text{once} \ \alp\in \big((k-1)\sqrt 3-2,k\sqrt 3-2),
\end{equation}
where $k=1,2,3,\cdots$ for $P=0,$ and $k=2,3$ for  $P>0.$
\item   If $\alp\in \big(2\sqrt 3-2,{3\over2}\sqrt{10}-2\big]\cup\big(3\sqrt 3-2,4\big)$ and $P=0$, then \eqref{R:2.1p4} still holds.
\item  If $\alp\in \big(2\sqrt 3-2,{3\over2}\sqrt{10}-2\big]\cup\big(3\sqrt 3-2,4\big)$ and $P>0$,
then
\begin{equation}\label{R:2.1p} \|V(t,\cdot )\|_{C^0(S^1)}\le MPe^{-(2-\mu ) t}.
\end{equation}
\end{enumerate}
\end{lem}

\noindent{\bf Proof.} Using Fourier analysis, we denote the Fourier series of $V$ and $\frac{\lam (3\m V^2+2V^3)}{\m ^3(V+\m )^2}$ as follows:

\begin{equation}\label{R:2.25}\arraycolsep=1.5pt\begin{array}{lll}
\hfill V(t,\theta )&=(2\pi)^{-\frac{1}{2}} \sum _{k\in \mathbb{Z}}a_k(t)e^{ik\theta},&\\[3mm]
\hfill \displaystyle\frac{\lam (3\m V^2+2V^3)}{\m ^3(V+\m )^2}&=(2\pi)^{-\frac{1}{2}}\sum _{k\in \mathbb{Z}}A_k(t)e^{ik\theta},&
\end{array}
\end{equation}
where $m>0$ is defined as in (\ref{R:mu}). It then follows from Lemma \ref{R:C} that $w(t,\theta )=e^{\varepsilon t}V(t,\theta)$ is bounded uniformly in $[t_0,+\infty)\times S^1$,
where $\varepsilon >0$ is the same as that of Lemma \ref{R:C}.
We thus obtain from (\ref{R:3.17'}) that $w$ satisfies the following evolution elliptic equation

\begin{equation}\label{R:2.27}
-w_{tt}+2(\mu +\varepsilon )w_t=w_{\theta\theta}-Pe^{-(2 -\mu -\varepsilon )t}
 +\Big[\mu ^2+\varepsilon ^2+2\mu \varepsilon +\displaystyle\frac{\lam (2\m +we^{-\varepsilon t})}{\m ^2(\m +we^{-\varepsilon t})^2}\Big]w
\end{equation}
in $(t_0,\infty)\times S^1.$
Similar to Lemma \ref{lemR.1}, one can derive from {\it a priori} estimates that $w$ and its derivatives, up to the third order, remain bounded on $[t_0,+\infty)\times S^1$, i.e.,

\begin{equation}\label{R:2.28}
\|e^{\varepsilon t}V(t,\theta  )\|_{C^3 ([t_0,+\infty)\times S^1 )}<+\infty.
\end{equation}
We thus obtain that there exists $C>0$ such that

\begin{equation}\label{R:2.29}
\sum _{k\in \mathbb{Z}}(k^2+1)|a_k(t)|^2  \le   Ce^{-2\varepsilon t},\qquad
\sum _{k\in \mathbb{Z}}(k^2+1)|A_k(t)|^2  \le   Ce^{-4\varepsilon t}.
\end{equation}
In view of (\ref{R:mu}), equation (\ref{R:3.17'}) can be rewritten as

\begin{equation}\label{R:3.17''}
-V_{tt}+ 2\mu V_t=V_{\theta\theta}+ 3\mu ^2V -\displaystyle\frac{\lam (3\m V^2+2V^3)}{\m ^3(V+\m )^2}-Pe^{-(2-\mu )t}.
\end{equation}
This implies  that $a_k(t)$ is a bounded solution of

\begin{equation}\label{R:2.31}
a_k''(t)-2\mu a'_k(t)+(3\mu ^2-k^2)a_k(t)=g_k(t):=\arraycolsep=1.5pt\left\{\begin{array}{lll}
A_0(t)+ \sqrt{2\pi}Pe^{-(2-\mu )t},\quad &k=0, \\[3mm]
A_k(t),\quad &k\neq0,
\end{array}\right.
\end{equation}
where $g_k(t)$ satisfies

\begin{equation}\label{R:2.31'}
|g_k(t)|\le\arraycolsep=1.5pt\left\{\begin{array}{lll}
M_0 e^{-2\varepsilon  t}+ \sqrt{2\pi}Pe^{-(2-\mu )t} ,\quad &k=0, \\[3mm]
\displaystyle\frac{M_0}{\sqrt {k^2+1}} e^{-2\varepsilon  t} ,\quad &k\neq0
\end{array}\right.
\end{equation}
for some constant $M_0>0$. 

Denote $d_{k,\mu}=k^2-2\mu ^2.$ By applying (\ref{R:2.29}), the integration of (\ref{R:2.31}) yields that

\[\arraycolsep=1.5pt\begin{array}{lll}
\hfill a_k(t)&=&\displaystyle\frac{1}{\sqrt{-d_{k,\mu}}}\bigg(e^{\mu t}\cos \sqrt{-d_{k,\mu}}t\displaystyle\int ^{+\infty}_te^{-\mu s}g_k(s)\sin \sqrt{-d_{k,\mu}}s\, ds
\\[2mm]
&&-e^{\mu t}\sin \sqrt{-d_{k,\mu}}t\displaystyle\int ^{+\infty}_te^{-\mu s} g_k(s)\cos \sqrt{-d_{k,\mu}}s\,ds\bigg),\quad \mbox{for} \; |k|<\sqrt 2\mu , \\[3mm]
a_k(t)&=&e^{\mu t} \displaystyle\int ^{+\infty}_t(s-t)e^{-\mu s}  g_k(s)ds,\quad \mbox{for} \; |k|=\sqrt 2 \mu,\\[3mm]
a_k(t)&=&\displaystyle\frac{1}{2\sqrt{d_{k,\mu}}}e^{(\mu -\sqrt {d_{k,\mu}})t} \displaystyle\int ^{+\infty}_t g_k(s)e^{-(\mu -\sqrt {d_{k,\mu}})s}ds\\
&&-\displaystyle\frac{1}{2\sqrt{d_{k,\mu}}}e^{(\mu +\sqrt {d_{k,\mu}})t} \displaystyle\int ^{+\infty}_t g_k(s)e^{-(\mu +\sqrt {d_{k,\mu}})s}ds,\quad \mbox{for} \; \sqrt 2\mu <|k|\le \sqrt 3\mu ,\\[3mm]
a_k(t)&=&\displaystyle a_k(t_0)e^{(\mu -\sqrt {d_{k,\mu}}) (t-t_0)}+\frac{e^{(\mu +\sqrt {d_{k,\mu}}) t_0}}{2\sqrt {d_{k,\mu}}}e^{(\mu -\sqrt {d_{k,\mu}}) (t-t_0)}\displaystyle\int_{t_0}^{+\infty}g_k(s)e^{-(\mu +\sqrt {d_{k,\mu}}) s}ds\\[3mm]
&&\displaystyle-\frac{1}{2\sqrt {d_{k,\mu}}}\displaystyle e^{(\mu -\sqrt {d_{k,\mu}}) t}\int_{t_0}^tg_k(s)e^{-(\mu -\sqrt {d_{k,\mu}}) s}ds\\[3mm]
&&\displaystyle-\frac{1}{2\sqrt {d_{k,\mu}}}\displaystyle e^{(\mu +\sqrt {d_{k,\mu}}) t}\int_t^{+\infty}g_k(s)e^{-(\mu +\sqrt {d_{k,\mu}}) s}ds,\quad \mbox{for} \; |k|> \sqrt 3\mu.
\end{array}\]
It then follows from the above that there exists $M_1>0$, depending only on $\alp$ and $P$, such that

\begin{equation}\label{FC}\arraycolsep=1.5pt\begin{array}{lll}
\hfill |a_0(t)|&\le & M_1\big(Pe^{-(2-\mu ) t}+e^{-2\varepsilon t}\big),  \\[3mm]
\hfill |a_k(t)|&\le & M_1e^{-2\varepsilon t}  \quad \mbox{for}\,\ 0<|k|\le \sqrt 3\mu, \\[3mm]
\hfill |a_{k}(t)|&\le & |a_k(t_0)|\, e^{-(\sqrt {d_{k,\mu}}-\mu ) (t-t_0)}
+\displaystyle\frac{M_1}{|k|\sqrt{k^2+1}}\Big(e^{-(\sqrt {d_{k,\mu}}-\mu )(t-t_0)}+e^{-2\varepsilon t}\Big),\quad \mbox{for}\quad |k|> \sqrt 3\mu.
\end{array}\end{equation}

In the following, we present proofs only for two special cases, as other cases can be proved in a similar way.

\vskip 0.05truein

{\bf  Case 1:\ \ $ \alp\in(\sqrt{3}-2,2\sqrt{3}-2).$}
In this case, we have

\begin{equation}\label{K1:dd}
1<\sqrt 3\mu<2,\quad    \sqrt {2^2-2\mu ^2}-\mu <2-\mu.
\end{equation}
It then follows from \eqref{R:2.29} and (\ref{FC}) that there exist positive constants $C_4,C_5,$ and $C_6$ such that

\begin{equation}\arraycolsep=1.5pt\begin{array}{lll}
\hfill &\|V(t, \cdot )\|^2_{H^1(S^1)}=\sum _{k\in \mathbb{Z}}(k^2+1)|a_k(t)|^2 \\[3mm]
\hfill =&|a_0(t)|^2 +\sum _{0<|k|<2}(k^2+1)|a_k(t)|^2+\sum _{|k|\geq 2}(k^2+1)|a_k(t)|^2 \\[3mm]
\hfill\leq & C_4 P^2 e^{-2(2-\mu)t}+C_5e^{-4\varepsilon t}+\sum _{|k|\geq 2}(k^2+1)C|a_k(t_0)|^2 e^{-2(\sqrt {2^2-2\mu ^2}-\mu ) (t-t_0)}\\[3mm]
&+\sum _{|k|\geq 2}(k^2+1)\displaystyle\frac{C}{k^2(k^2+1)}\Big(e^{-2(\sqrt {2^2-2\mu ^2}-\mu ) (t-t_0)}+e^{-4\varepsilon t}\Big)\\[3mm]
\leq & C_4 P^2e^{-2(2-\mu)t}+C_5 e^{-4\varepsilon t}+C_6 e^{-2(\sqrt {2^2-2\mu ^2}-\mu ) t}.
\end{array}\label{4:4.21}\end{equation}
By (\ref{K1:dd}), this further implies that both for $P=0$ and $P>0,$

\begin{equation}\label{K1:d}\|V(t, \cdot )\|_{C^0(S^1)}\le N_1 e^{-2\varepsilon t}+N_2 e^{-(\sqrt {2^2-2\mu ^2}-\mu ) t}, \end{equation}
 where $N_1$ and $N_2$ are positive constants.
If $\sqrt {2^2-2\mu ^2}-\mu \le 2\varepsilon  $, then (\ref{K1:d}) implies (\ref{R:2.1p4}), and the proof of Case 1 is thus complete. Otherwise, we repeat the above procedure with $\varepsilon$ replaced by $2 \varepsilon$ in (\ref{R:2.1}). By taking finite similar steps, we reach a finite integer $n$ such that $2^n\varepsilon \ge \sqrt {2^2-2\mu ^2}-\mu$, and the estimate (\ref{K1:d}) holds for $2 \varepsilon$ replaced by $2^n\varepsilon$. Therefore, we conclude that there exists $C>0$ such that

$$\|V(t, \cdot )\|_{C^0(S^1)}\le C e^{-(\sqrt {2^2-2\mu ^2}-\mu ) t},$$
and hence the estimate (\ref{R:2.1p4}) is proved in this case.

\vskip 0.05truein

{\bf  Case 2:\ \ $\alp\in \big(2\sqrt 3-2,3\sqrt 3-2\big).$} In this case, we have $2< \sqrt 3\mu <3$ and

\begin{equation}\label{530}
 2-\mu\left\{\arraycolsep=1.5pt
\begin {array}{lll}
<\sqrt {3^2-2\mu ^2}-\mu,\quad \text{for}\  \alp \in  (2\sqrt 3-2, {3\over2}\sqrt{10}-2), \\
=\sqrt {3^2-2\mu ^2}-\mu,\quad \text{for}\  \alp= {3\over2}\sqrt{10}-2,\\
> \sqrt {3^2-2\mu ^2}-\mu, \quad \text{for}\  \alp \in  ({3\over2}\sqrt{10}-2, 3\sqrt 3-2).
\end{array}\right.
\end{equation}
We then deduce from (\ref{FC}) that there exist constants $C_i>0$ ($i=9, 10, 11, 12$) such that

\begin{equation}\arraycolsep=1.5pt\begin{array}{lll}
\hfill &\|V(t, \cdot )\|^2_{H^1(S^1)}=\sum _{k\in \mathbb{Z}}(k^2+1)|a_k(t)|^2 \\[3mm]
\hfill =&|a_0(t)|^2 +\sum _{0<|k|<3}(k^2+1)|a_k(t)|^2+\sum _{|k|\geq 3}(k^2+1)|a_k(t)|^2 \\[3mm]
\leq & C_9 P^2e^{-2(2-\mu)t}+C_{10} e^{-4\varepsilon t}+C_{11} e^{-2(\sqrt {3^2-2\mu ^2}-\mu ) t}.
\end{array}\end{equation}
By \eqref{530}, this further implies that

\begin{equation}\label{K1:c}
\|V(t, \cdot )\|_{C^0(S^1)}\le\left\{\arraycolsep=1.5pt
\begin {array}{lll}

N_3e^{-2\varepsilon t}+N_4 e^{-(\sqrt {3^2-2\mu ^2}-\mu ) t},\quad &\text{for}\ \alp\in  (2\sqrt 3-2, {3\over2}\sqrt{10}-2] \ \text{and}\ P=0,\\& \text{or}\ \ \alp\in [{3\over2}\sqrt{10}-2, 3\sqrt 3-2)\ \text{and}\ P\geq 0, \\
 N_5e^{-2\varepsilon t}+N_6 Pe^{-(2-\mu ) t},\quad &\text{for}\ \alp\in  (2\sqrt 3-2, {3\over2}\sqrt{10}-2] \ \text{and}\ \ P>0,
 \end{array}\right.
 \end{equation}
where $N_i, i=3,4,5,6$ are positive constants. Furthermore, similar to Case 1, one can obtain the estimate \eqref{R:2.1p4} for $P=0$ and (\ref{R:2.1p}) for $P>0$.
This completes the proof of Lemma \ref{R:D}.
\qed

\begin{prop} \label{prop:K1} Under the assumptions of Theorem \ref{thm3}, suppose $V(t,\theta)$ satisfies (\ref{R:newtrans}) and $\mu $ is given by (\ref{R:mu}).
Then we have the following results:
\begin{enumerate}

\item  If $\alp\in \mathring{\mathcal{A} }\backslash\big\{(2\sqrt{3}-2,{3\over2}\sqrt{10}-2]\cup(3\sqrt{3}-2,4)\big\},$ then for both $P=0$ and $P>0,$ once $\alp\in \big((k-1)\sqrt 3-2,k\sqrt 3-2\big), $ there exist $A_k\in \R $ and $\theta _k\in S^1$ such that
\begin{equation}\label{4.3:1a}
\lim_{t\to+\infty} e^{(\sqrt {k^2-2\mu ^2}-\mu) t}V(t,\theta ) = A_k\sin (k\theta +\theta _k)\quad \mbox{in}\ C^2(S^1),
\end{equation}
where $k=1,2,3,4,\cdots$ for $P=0,$ and $k=1,2,3$ for  $P>0.$
\item   If $\alp\in \big(2\sqrt 3-2,{3\over2}\sqrt{10}-2\big]\cup\big(3\sqrt 3-2,4\big)$ and $P=0$, then \eqref{4.3:1a} still holds.

\item If $\alp\in \big(2\sqrt 3-2,{3\over2}\sqrt{10}-2\big)\cup\big(3\sqrt 3-2,4\big)$ and $P>0$,
then

\begin{equation}\label{4.3:1}
\lim_{t\to+\infty}e^{(2-\mu ) t} V(t,\theta )=\frac{9P}{36+2(2+\alp )^2}\quad \mbox{in}\ C^2(S^1).
\end{equation}
\item If $\alp={3\over2}\sqrt{10}-2$ and $P>0$, then there exist $A_3\in \R $ and $\theta _3\in S^1$ such that

\begin{equation}\label{76}
\lim_{t\to+\infty}e^{(\sqrt {3^2-2\mu ^2}-\mu) t}V(t,\theta )=A_3\bigg(\sin (3\theta +\theta _3)+{P\over 9}\bigg)\quad \mbox{in}\ C^2(S^1).
\end{equation}
\end{enumerate}

\end{prop}

\noindent{\bf Proof.}  Define

\begin{equation}\label{G:2.39}
w(t,\theta )=e^{\gamma t}V(t,\theta),
\end{equation}
where  $\gamma >0$ is to be chosen later, such that $w$ satisfies

\begin{equation}\label{G:2.40}
-w_{tt}+2(\gamma +\mu )w_t=w_{\theta\theta} -Pe^{[\gamma -(2-\mu)]t}
+\Big[\mu ^2+\gamma (\gamma +2\mu )+\displaystyle\frac{\lam (2\m +we^{-\gamma t})}{\m ^2(\m +we^{-\gamma t})^2}\Big]w
\end{equation}
in $(t_0,+\infty)\times S^1.$ Note that

\beq
\frac{\lam (2\m +we^{-\gamma t})}{\m ^2(\m +we^{-\gamma t})^2}-2\mu^2=O(e^{-\gamma t})\quad \text{as}\ t\rightarrow +\infty.
\eeq
In the following, we only provide proofs for two special cases, as other cases can be proved in a similar way.
\vskip 0.05truein

{\bf  Case 1:\ \ $ \alp\in(\sqrt{3}-2,2\sqrt{3}-2).$}
In this case, we take $\gamma =\sqrt {2^2-2\mu ^2}-\mu <2-\mu.$ Then $\mu ^2+\gamma (\gamma +2\mu )=4-2\mu^2$ and $w$ is bounded uniformly in $[t_0,+\infty)\times S^1$ in view of Lemma \ref{R:D}. Moreover, it follows from (\ref{G:2.40}) that

\begin{equation}\label{G:2.40a}
-w_{tt}+2\sqrt{2^2-2\mu ^2} w_t
=w_{\theta\theta}
 +4w-e^{-\gamma t}f(t, \theta )-Pe^{-[2-\sqrt {2^2-2\mu ^2}]t}\quad \mbox{in}\ \ (0,+\infty)\times S^1.
\end{equation}
Similar to Lemma \ref{lemR.1},
if we
define the ``$\omega$-limit set" $\mathrm{\Gamma}(\mathcal{L}'')$ of the ``orbit" $\mathcal{L}'':=\{w(t,
\cdot ):\, t\ge t_0\}$ for (\ref{G:2.40}) with $\gamma =\sqrt {2^2-2\mu ^2}-\mu$ as

\[
\mathrm{\Gamma}(\mathcal{L}''):=  \bigcap _{t\ge t_0}\overline{\big\{w(\iota ,\cdot );\, \iota \ge t\big\}} ,
\]
where the closure is with respect to the topology of $C^2(S^1)$, then we obtain that $\mathrm{\Gamma}(\mathcal{L}'')$ is a nonempty, compact,
and connected set in $C^2(S^1)$. Moreover, $\mathrm{\Gamma}(\mathcal{L}'') \subset \mathfrak{S}''$, where
$\mathfrak{S}''$ is a nonempty, compact, and connected subset of

\begin{equation}\label{73}
\Big\{\psi (\theta ) \in C^2(S^1):\, \frac{d^2\psi}{d\theta ^2}+4\psi =0\Big\}=\Big\{A_0\sin (2\theta +\theta _0): \, A_0\in \mathbb{R},\ \ \theta _0\in S^1 \Big\}.
\end{equation}

We next further analyze the limit behavior of $w(t,\cdot)$ as $t\to\infty$.
Consider the bounded Fourier coefficients $a_{k}(t)$ of $w(t,\theta )$, which are defined by

\[
a_{k}(t):=(2\pi )^{-1/2}\int _{S^1}w(t,\theta )e^{ - ik \theta}d\theta , \quad k\in \mathbb{Z}.
\]
It then follows from (\ref{G:2.40a}) that

\begin{equation}\label{G:2.42}
a_k''(t)-2\sqrt{2^2-2\mu ^2} a_k'(t)+(4 -k^2)a_k(t)=F_k(t):=\arraycolsep=1.5pt\left\{\begin{array}{lll}
e^{-\gamma t}f_0(t)+ \sqrt{2\pi}Pe^{-[2-\sqrt {2^2-2\mu ^2}]t},\quad &k=0, \\[3mm]
e^{-\gamma t}f_k(t),\quad &k\neq0,
\end{array}\right.
\end{equation}
  where  $f_k(t)=(2\pi )^{-1/2}\int _{S^1}f(t,\theta)e^{ - ik \theta}d\theta$ is  bounded  uniformly in $[t_0,+\infty)$. Also denote $d_{k,\mu}=k^2-2\mu ^2$ as in Lemma \ref{R:D}.
By the uniform boundedness of $a_k(t)$, the integration of (\ref{G:2.42}) yields that

\begin{eqnarray*}
\hfill a_k(t)&=&\displaystyle\frac{1}{\sqrt{-d_{k,\mu}}}\bigg(e^{\sqrt{d_{2,\mu}} t}\cos \sqrt{-d_{k,\mu}}t\displaystyle\int _t^{+\infty}e^{-\sqrt{d_{2,\mu}} s} F_k(s)\sin \sqrt{-d_{k,\mu}}s\,ds\\[2mm]
&&-e^{\sqrt{d_{2,\mu}} t}\sin \sqrt{-d_{k,\mu}}t\displaystyle\int _t^{+\infty}e^{-\sqrt{d_{2,\mu}} s}F_k(s)\cos \sqrt{-d_{k,\mu}}s\, ds\bigg),\quad \mbox{for}\;\; |k|<\sqrt 2\mu , \\[3mm]
\quad a_k(t)&=&e^{\sqrt{d_{2,\mu}} t} \displaystyle\int _t^{+\infty}(s-t)e^{-\sqrt{d_{2,\mu}} s}  F_k(s)ds,\quad \mbox{for}\;\; |k|=\sqrt 2 \mu,\\[3mm]
\quad a_k(t)&=&\displaystyle\frac{1}{2\sqrt{d_{k,\mu}}}e^{(\sqrt{d_{2,\mu}} -\sqrt {d_{k,\mu}})t}
 \displaystyle\int
_t^{+\infty}\Big[e^{2\sqrt {d_{k,\mu}}s}-e^{2\sqrt {d_{k,\mu}}t}\Big] \\[3mm]
&&\qquad \cdot e^{-(\sqrt{d_{2,\mu}} +\sqrt {d_{k,\mu}})s}F_k(s)ds,\quad \mbox{for}\;\; \sqrt 2\mu <|k|<   2 ,\\[3mm]
 a_{\pm2}(t)&=&B_{\pm2}+\displaystyle\frac{1}{2\sqrt{d_{2,\mu}}}\displaystyle\int
_t^{+\infty}\bigg(1-e^{2\sqrt{d_{2,\mu}}(t-s)}\bigg)
F_{\pm2}(s)ds,\quad \mbox{for}\;\;
|k|=  2 ,\\[3.5mm]
\hfill  a_k(t)&=&\displaystyle a_k(t_0)e^{(\sqrt{d_{2,\mu}} -\sqrt {d_{k,\mu}}) (t-t_0)}
+e^{(\sqrt{d_{2,\mu}} -\sqrt {d_{k,\mu}}) (t-t_0)}\\[3mm]
&&
\cdot\displaystyle\frac{e^{(\sqrt{d_{2,\mu}}+\sqrt {d_{k,\mu}})  t_0}}{2\sqrt {d_{k,\mu}}
}\int_{t_0}^{+\infty}e^{-(\sqrt{d_{2,\mu}}+\sqrt {d_{k,\mu}}) s}F_k(s)ds\\[3mm]
&&\displaystyle-\frac{e^{(\sqrt{d_{2,\mu}} -\sqrt {d_{k,\mu}}) t}}{2\sqrt {d_{k,\mu}}}\int_{t_0}^te^{-(\sqrt{d_{2,\mu}} -\sqrt {d_{k,\mu}}) s}F_k(s)ds\\[3.5mm]
&&\displaystyle-\frac{e^{(\sqrt{d_{2,\mu}}+\sqrt {d_{k,\mu}}) t}}{2\sqrt {d_{k,\mu}}
}\int_t^{+\infty}e^{-(\sqrt{d_{2,\mu}}+\sqrt {d_{k,\mu}})s}F_k(s)ds,\quad \mbox{for}\;\; |k|> 2,\\
\end{eqnarray*}
where $B_{\pm2}$ are two complex constants satisfying

\beq
B_{\pm2}=a_{\pm2}(t_0)-\displaystyle\frac{1}{2\sqrt{d_{2,\mu}}}
\cdot \displaystyle\int
_{t_0}^{+\infty}\bigg(1-e^{2\sqrt{d_{2,\mu}}(t_0-s)}\bigg)
F_{\pm2}(s)ds.
\eeq
We calculate from the above that there exists a constant $C_{13}>0$ such that

\[\arraycolsep=1.5pt\begin{split}
 |a_0(t)| &\le   C_{13}e^{- \beta t},\, \ \mbox{where}\,\ \beta=\min\big\{\sqrt {d_{2,\mu}}-\mu,\, 2-\sqrt {d_{2,\mu}}\big\},  \\
 |a_k(t)| &\le  C_{13}e^{- \gamma t}, \quad   |k|<2, \\
 |a_{\pm2}(t)-B_{\pm2}| &\le C_{13}e^{- \gamma t}, \\
 |a_k(t)|  &\le  C_{13}\Big[e^{-(\sqrt {d_{k,\mu}}-\sqrt {d_{2,\mu}})t}+e^{- \gamma t}\Big], \quad   |k|>2.
\end{split}
\]
Therefore, $a_k(t)\rightarrow 0$ exponentially as $t\rightarrow +\infty$ if $|k|\neq 2$, and $a_{\pm 2}(t)\rightarrow a  _{\pm }\in\R $  exponentially as $t\rightarrow +\infty$.  We then conclude from (\ref{73}) and the above that there exist
$A_2\in \R$ and $\theta _2\in S^1$ such that

\begin{equation}
 \lim_{t\to+\infty} e^{(\sqrt {2^2-2\mu ^2}-\mu) t} V(t,\theta )=A_2\sin (2\theta +\theta _2) \quad \mbox{as}\quad t\to+\infty ,\label{sing:1B}
\end{equation}
which is a special case of (\ref{4.3:1a}).

\vskip 0.05truein

\vskip 0.05truein

{\bf  Case 2:\ \  $\alp \in  (2\sqrt 3-2, {3\over2}\sqrt{10}-2]$ and $P>0$.}
In this case, we take $\gamma =2-\mu$, such that $w$ is bounded uniformly in $[t_0,+\infty)\times S^1$ in view of Lemma \ref{R:D}.
We then derive from (\ref{G:2.40}) that

\begin{equation}\label{G:2.40B}
-w_{tt}+4w_t=w_{\theta\theta}-P+(4+2\mu^2)w -e^{-(2-\mu) t}f(t, \theta )\quad \mbox{in}\ \ (0,+\infty)\times S^1,
\end{equation}
where $f(t,\theta)$ is bounded uniformly in $[0,+\infty)\times S^1$. Moreover, we know that in this case, with $\gamma =2-\mu$,
$\mathrm{\Gamma}(\mathcal{L}'') \subset \mathfrak{S}''$, where
$\mathfrak{S}''$ is a nonempty, compact, and connected subset of

\begin{equation}\label{integer:1}
\begin{split}
&\Big\{\psi (\theta ) \in C^2(S^1):\, \frac{d^2\psi}{d\theta ^2}+(4+2\mu^2)\psi-P =0\Big\}\\=&\Big\{A_0\sin (\sqrt{4+2\mu^2}\theta +\theta _0)+{P\over{4+2\mu^2}}: \, A_0\in \mathbb{R},\ \ \theta _0\in S^1 \Big\}\\
=&
\left\{\arraycolsep=1.5pt
\begin {array}{lll}
\Big\{{P\over{4+2\mu^2}}\Big\},\quad&\text{for}\ \alp \in  (2\sqrt 3-2, {3\over2}\sqrt{10}-2),\\[1.5mm]
\Big\{A_0\sin (3\theta +\theta _0)+{P\over 9}: \, A_0\in \mathbb{R},\ \ \theta _0\in S^1 \Big\}, \quad&\text{for}\ \alp={3\over2}\sqrt{10}-2,
\end{array}\right.
\end{split}
\end{equation}
because $\sqrt{4+2\mu^2}$ is not an integer for $\alp \in  (2\sqrt 3-2, {3\over2}\sqrt{10}-2),$ and $\sqrt{4+2\mu^2}=3$ for $\alp={3\over2}\sqrt{10}-2$. Therefore, similar to Case 1, one can further derive from (\ref{integer:1}) that $V(t,\theta )$ satisfies \eqref{4.3:1} for $\alp \in  (2\sqrt 3-2, {3\over2}\sqrt{10}-2),$ and satisfies \eqref{76} for $\alp={3\over2}\sqrt{10}-2.$
This completes the proof of Proposition \ref{prop:K1}.
\qed

We finally remark that Theorem  \ref{thm3} follows immediately from Proposition \ref{prop:K1} and (\ref{R:newtrans}).

\appendix

\renewcommand{\appendixname}{Appendix}

\begin{appendices}
\section{}
In this appendix, inspired by \cite{GZ,HW} we complete the proof of Lemma \ref{thm3.1}. We first establish the following estimates.

\begin{lem}\label{23}
Suppose $\phi$ is a nonnegative smooth function satisfying
$$-\Delta\phi+\frac{2(2+\alpha)^2}{9}\phi_{x_1}\leq \lambda \phi^4+P\phi^2\,\ \hbox{in}\,\ \overline{B_R(x_0)}\subset \mathbb{R}^2,$$
where $x_0\in\R^2$, $0<R\leq 1$, $\alpha >-2$, $\lambda>0$, and $P\ge 0$ are given constants. Then there exists a constant $\eta_0>0$, depending only on $\alpha$, $\lambda$, and $P$, such that the estimate

\begin{equation}\label{22}
\frac{1}{r^{4/3}}\int_{B_r(y)}\phi dx\leq \eta_0, \,\ \forall B_r(y)\subset B_R(x_0),
\end{equation}
implies that

\begin{equation} \label{22A}
\phi(x)\leq\frac{2}{R} \,\ \hbox{in}\,\ B_{\frac{R}{2}}(x_0).
\end{equation}
\end{lem}

\noindent{\bf Proof.}
Inspired by \cite[Lemma 2.1]{GZ} and \cite[Lemma 2.2]{HW}, we denote $K=\max_{|x-x_0|\leq R}(R-|x-x_0|)\phi(x)>0$, where $0<R\le 1$. Choose $\xi\in B_{R}(x_0)$ such that $(R-|\xi-x_0|)\phi(\xi)=K$, and set $\sigma=R-|\xi-x_0|.$ Then we have $\phi(x)\leq \frac{2K}{\sigma}$ for $x\in B_{\frac{\sigma}{2}}(\xi)$.  Denote $\mu=\frac{K}{\sigma}=\phi(\xi)$ and consider $\psi(x)=\frac{1}{\mu}\phi\big(\xi+\frac{1}{\mu^{3/2}}x\big)$, such that $\psi(x)$ satisfies

\begin{equation}\label{21}
\begin{cases}
-\Delta \psi+\frac{1}{\mu^{3/2}}\frac{2(2+\alpha)^2}{9}\psi_{x_1}\leq \lambda \psi^4 +\frac{P}{\mu^2} \psi^2 \ \ \text{in}\,\ B_{\sigma\mu^{3/2}}(0),\\
0\leq\psi\leq2 \quad \text{in}\,\ B_{\frac{1}{2}\sigma\mu^{3/2}}(0),\\
\psi(0)=1.
\end{cases}
\end{equation}
One can note that if $K\leq 1$ holds, then (\ref{22A}) follows immediately.

It now suffices to prove that $K\leq 1$ holds. First, if $\mu\leq 1,$ then it is clear that $K=(R-|\xi-x_0|)\mu\leq R\leq1.$ It only remains therefore to prove that $K\leq 1$ for the case where $\mu>1.$
On the contrary, suppose $\mu>1$ and $K>1.$  We then have $\sigma\mu^{3/2}=K\mu^{1/2}>1$ and it thus follows from \eqref{21} that

\begin{equation}
\begin{cases}
-\Delta \psi+\frac{1}{\mu^{3/2}}\frac{2(2+\alpha)^2}{9}\psi_{x_1}\leq (8\lambda +2P) \psi  \quad \text{in}\ \ \overline{B_{\frac{1}{2}}(0)},\\
0\leq\psi\leq2 \quad \text{in}\ B_{\frac{1}{2}}(0),\\
\psi(0)=1,\,\ \frac{1}{\mu^{3/2}}\in (0,1)\,\ \text{is bounded}.
\end{cases}\label{3:3.7}
\end{equation}
Moreover, we obtain from \eqref{22} that

\[
\int_{B_{\frac{1}{2}}(0)}\psi(x)dx=\mu^2\int_{B_{\frac{1}{2\mu^{3/2}}}(\xi)}\phi(y)dy\leq 2^{-\frac{4}{3}}\eta_0,
\]
because $\frac{1}{2\mu^{3/2}}<\frac{1}{2\mu}=\frac{\sigma}{2K}<\frac{\sigma}{2}$ for $\mu>1$ and $K>1.$
By the elliptic estimate \cite[p. 244]{GT}, we then deduce from (\ref{3:3.7}) that

 \begin{equation}\label{3:3.9}
\psi(0)\leq C \int_{B_{\frac{1}{2}}(0)}\psi(x)dx\leq 2^{-\frac{4}{3}}C\eta_0,
\end{equation}
where $C$ is a constant depending only on $\alpha$, $\lambda$, and $P$. By choosing $\eta_0>0$ small enough that $2^{-\frac{4}{3}}C\eta_0<1,$ we conclude from (\ref{3:3.9}) that $\psi(0)<1$, which is a contradiction in view of (\ref{21}). This shows that $K\leq 1$ also holds for the case where $\mu>1$, and we have finished.
\qed

\noindent{\bf Proof of Lemma \ref{thm3.1}.} Under the assumption (\ref{Z:2}), we first prove that there exist  $C>0$ and $t_1\geq 0$ such that

\begin{equation} \label{3:upbound}
v(t,\theta)\geq C, \,\ \forall (t,\theta)\in(t_1,+\infty)\times S^1,
\end{equation}
where $C$ depends only on $\alpha,\lambda,\beta$, and $C_\beta$ given in \eqref{Z:2}. By taking $w=\frac{1}{v}$, it follows from \eqref{R:trans-1} that

\[
-(w_{tt}+w_{\theta\theta})+\frac{2(2+\alpha)^2}{9}w_t=\lambda w^4+Pe^{-\frac{4-\alpha}{3}t}w^2-\frac{(2+\alpha)^2}{9}w-2\frac{w_t^2+w_\theta^2}{w}\leq \lambda w^4+Pw^2.
\]
Let $\eta_0>0$ be the same as that of Lemma \ref{23}. Under the assumption \eqref{Z:2}, we then have

\[
\frac{1}{r^{4/3}}\int_{B_r(x)}w(y)dy\leq C_\beta r^\beta\leq\eta_0,\quad \forall r<\bar{r},\,\ x=(t,\theta)\in(t_0,+\infty)\times S^1,
\]
where the constant $\bar{r}=\min\{\big(\frac{\eta_0}{C_\beta}\big)^\frac{1}{\beta},1\}>0$ depends only on $\beta, \alpha ,\lambda$, and $P$. Applying Lemma \ref{23}, we then derive that there exists $t_1>0$ such that

\[
w(t,\theta)\leq \frac{2}{\bar{r}} ,  \,\ \forall (t,\theta)\in(t_1,+\infty)\times S^1,
\]
from which the estimate (\ref{3:upbound}) follows.

To complete the proof of Lemma \ref{thm3.1}, we only need to prove that $u(x)\leq C|x|^{\frac{2+\alpha}{3}}$ near the origin. By applying (\ref{3:upbound}), it is standard to derive (e.g., \cite[Lemma 2.3]{GZ}) that
there exists a constant $C>0$, depending only on $\alpha,\lambda, P,\beta$, and $C_\beta,$ such that the following spherical Harnack inequality holds

\begin{equation} \label{3:lowbound}
\sup_{|x|=r}u(x)\leq C\inf_{|x|=r}u(x),\quad \forall r\in\big(0,\frac{1}{2}\big).
\end{equation}
Define

\[
\bar{u}(r)=\frac{1}{2\pi}\int_0^{2\pi}u(r,\theta)d\theta,
\]
where $u(r,\theta)=u(x)$, such that by \eqref{62}, $\bar{u}$ satisfies

\[
(r\bar{u}_r)_r=\frac{r}{2\pi}\int_0^{2\pi}\Big(\frac{\lambda r^\alpha}{u^2}+P\Big)d\theta=\frac{r}{2\pi}\int_0^{2\pi}\Big(\frac{\lambda r^{\alpha-\frac{2(2+\alpha)}{3}}}{v^2}+P\Big)d\theta=\frac{1}{2\pi}\int_0^{2\pi}\Big(\frac{\lambda r^{\frac{\alpha-1}{3}}}{v^2}+Pr\Big)d\theta.
\]
This implies that $r\bar{u}_r$ is increasing monotonically near the origin, and thus the limit $\lim_{r\rightarrow0^+}r\bar{u}_r$ exists. Moreover,
 it cannot be negative because otherwise $\bar{u}$ must be negative near the origin, which is impossible.
  We now prove that $\lim_{r\rightarrow0^+}r\bar{u}_r=0.$ If it is false, then $r\bar{u}_r\geq C>0,$ i.e., $\bar{u}_r>\frac{C}{r},$ near the origin. This implies that

\[
+\infty>\bar{u}(1)=\bar{u}(1)-\bar{u}(0)=\int_0^{1}\bar{u}_rdr>C\int_0^{1}\frac{1}{r}dr=+\infty,
\]
which is a contradiction. Therefore, we have $\lim_{r\rightarrow0^+}r\bar{u}_r=0.$

By estimate (\ref{3:upbound}), we have

\[
(r\bar{u}_r)_r=\frac{1}{2\pi}\int_0^{2\pi}\Big(\frac{\lambda r^{\frac{\alpha-1}{3}}}{v^2}+Pr\Big)d\theta\leq C\lambda r^{\frac{\alpha-1}{3}}+Pr.
\]
Given any $\varepsilon>0,$ integrating the above inequality from $\varepsilon$ to $r$, we obtain that

\[
r\bar{u}_r(r)-\varepsilon\bar{u}_r(\varepsilon)\leq \frac{3C\lambda}{\alpha+2}(r^{\frac{\alpha+2}{3}}-\varepsilon^{\frac{\alpha+2}{3}})+{\frac{P}{2}r^2-\frac{P}{2}\varepsilon^2}.
\]
Because $\lim_{r\rightarrow0^+}r\bar{u}_r=0$, the above estimate gives that

\[
r\bar{u}_r(r)\leq \frac{3C\lambda}{\alpha+2}r^{\frac{\alpha+2}{3}}+\frac{P}{2}r^2,
\]
i.e.,

\[
\bar{u}_r(r)\leq \frac{3C\lambda}{\alpha+2}r^{\frac{\alpha-1}{3}}+\frac{P}{2}r.
\]
We thus obtain that

\begin{equation}\label{27}
\bar{u}(r)=\bar{u}(r)-\bar{u}(0)=\int_0^r\bar{u}(r)_rdr\leq \frac{9C\lambda}{(\alpha+2)^2}r^{\frac{\alpha+2}{3}}+\frac{P}{4}r^2.
\end{equation}
By (\ref{3:lowbound}) and \eqref{27} we conclude that

\[
u(x)\leq\sup_{|y|=|x|}u(y)\leq C\inf_{|y|=|x|}u(y)\leq C\bar{u}(|x|)\leq \frac{9C\lambda}{(\alpha+2)^2}|x|^{\frac{\alpha+2}{3}}+\frac{P}{4}|x|^2=
\Big(\frac{9C\lambda}{(\alpha+2)^2}+\frac{P}{4}|x|^{\frac{4-\alpha}{3}}\Big)|x|^{\frac{\alpha+2}{3}}.
\]
Because either $\alp \geq0$ for $P=0$ or $0\leq \alp \le 4$ for $P>0$, the above estimate gives that $u(x)\leq C|x|^{\frac{2+\alpha}{3}}$ holds near the origin, and the proof is therefore complete.
\qed

\section{}

Recall that the set $\mathfrak{S}$, which is defined by \eqref{2.6S}, denotes the set of all positive solutions for \eqref{2.6}.
Define the following  functional

\beq\label{52}
E(v)=\int_{S^1}\Big(\frac{1}{2}v_{\theta}^2-\frac{A^2}{2}v^2-\frac{\lambda}{v}\Big)d\theta ,
\eeq
where $A>0$ and $\lam >0$ are as in Section 3. In this appendix, we derive the following $\L$ojasiewicz--Simon-type inequality in terms of $E(\cdot)$.

\begin{lem}\label{516}
For any $w\in \mathfrak{S}$, which is defined by \eqref{2.6S}, there exist positive constants $\sigma>0$ and $\bbtheta\in(0,{1\over2})$, depending only on $w$, such that for all $v\in H^2(S^1)$ and $\|v-w\|_{H^2(S^1)}<\sigma$,

\beq\label{519}
\big\|-v_{\theta\theta}-A^2 v+\frac{\lambda}{v^2}\big\|_{L^2(S^1)}\geq |E(v)-E(w)|^{1-\bbtheta},
\eeq
where $E(v)$ is defined by \eqref{52}.
\end{lem}

\noindent{\bf Proof.} Our proof is inspired by \cite{Z02}. We first consider the linearized problem of \eqref{6} near the equilibrium $w\in\mathfrak{S}$:

\[
L\varphi=-\varphi_{\theta\theta}-A^2\varphi-\frac{2\lambda}{w^3}\varphi,\,\ \varphi\in H^2(S^1).
\]
It is easy to see that the operator $L$ defined on $H^2(S^1)\subset L^2(S^1)$ is a self-adjoint operator. Define the bilinear form $B[ \cdot ,\cdot ]$ by

\[
B[h,k]=(L h,k)_{L^2(S^1)}=\int_{S^1}h_{\theta}k_{\theta}d\theta-A^2\int_{S^1}hkd\theta-
2\lambda\int_{S^1}\frac{1}{w^3}hkd\theta,\,\ h, k\in H^1(S^1).
\]
Because it follows from \eqref{4:1} that $0<C_1<w(\theta)<C_2$ on $S^1$, we have

\[
\big|B[h,k]\big|\leq \|h_{\theta}\|_{L^2(S^1)}\|k_{\theta}\|_{L^2(S^1)}+C \|h\|_{L^2(S^1)} \|k\|_{L^2(S^1)}\leq C \|h\|_{H^1(S^1)} \|k\|_{H^1(S^1)},
\]
and

\beq \label{51}
\|h_{\theta}\|^2_{L^2(S^1)}\leq B[h,h]+C\|h\|^2_{L^2(S^1)},
\eeq
which then implies that

\[
\|h\|^2_{H^1(S^1)}\leq B[h,h]+(C+1)\|h\|^2_{L^2(S^1)}.
\]
Thus, there exists a real constant $\gamma>0$ such that the operator $\gamma I+L$ is coercive on $H^1(S^1).$
Using the Lax--Milgram theorem, a Fredholm alternative result then holds for the problem

\[
L\varphi=h,\,\ \varphi\in H^2(S^1)\subset L^2(S^1).
\]
More precisely, we have either $ker(L)=\emptyset$ or $\dim(ker L)=m>0$ for some $m\in \mathbb{N}$, in which case the equation $L\varphi=h$ has a solution if and only if
$h\in(ker L)^\perp.$

We now focus on the case $\dim(ker L)=m>0$ to finish the proof of the lemma. Let $(\varphi_1,\varphi_2,\cdots,\varphi_m)$ be the normalized orthogonal basis of $ker(L)$ in $L^2(S^1)$, and denote by $\Pi$ the projection from $L^2(S^1)$ onto $ker(L)$. Define the operator $\mathscr{L}$ from $H^2(S^1)$ onto $L^2(S^1)$ as follows:

\[
\mathscr{L}\varphi=\Pi \varphi+L \varphi,\,\ \varphi\in H^2(S^1).
\]
Then

\[
\mathscr{L}:H^2(S^1)\mapsto L^2(S^1)
\]
is a one-to-one and onto operator.
Define $\psi=v-w$ and

\beq\label{55}
\mathscr{M}\psi=-v_{\theta\theta}-A^2 v+\frac{\lambda}{v^2}: H^2(S^1)\mapsto L^2(S^1).
\eeq
Note that $\mathscr{M}$ is Frechet differentiable at $\psi=0$ because $w>0$ on $S^1.$
It is easy to see that

\[
D\mathscr{M}(0)=L,
\]
where $D\mathscr{M}$ denotes the Frechet derivative of $\mathscr{M}.$ Denote

\[
\mathscr{N}\psi=\mathscr{M}\psi+\Pi \psi,\,\ \psi\in H^2(S^1),
\]
such that

\[
D\mathscr{N}(0)=\mathscr{L}.
\]

Because $\mathscr{L}$ is a one-to-one and onto operator, by the local inversion theorem (see, e.g. \cite[\S2.7]{N}) in nonlinear analysis, there exist a small neighborhood $W_1(0)$ of the origin in $H^2(S^1)$, a small neighborhood $W_2(0)$ of the origin in $L^2(S^1)$, and an inverse mapping

\[
\mathscr{T}:\, W_2(0)\mapsto W_1(0),
\]
such that

\[
\mathscr{N}(\mathscr{T}(g))=g,\quad \forall g\in W_2(0),
\]
and

\[
\mathscr{T}(\mathscr{N}(\psi))=\psi,\quad \forall \psi\in W_1(0).
\]
Because $0<C_1<w$ on $S^1$, $\frac{1}{v^2}=\frac{1}{(w+\psi)^2}$ is analytic in $\psi\in W_1(0)$. Thus, the operator  $ \mathscr{N}$ and its inverse mapping $\mathscr{T}$ are all analytic. Furthermore, there exists a positive constant $C>0$ such that

\beq \label{57}
\|\mathscr{T}(g_1)-\mathscr{T}(g_2)\|_{H^2(S^1)}\leq C \|g_1-g_2\|_{L^2(S^1)},\quad \forall g_1,g_2\in W_2(0),
\eeq
and

\beq\label{56}
\|\mathscr{N}(\psi_1)-\mathscr{N}(\psi_2)\|_{L^2(S^1)}\leq C\|\psi_1-\psi_2\|_{H^2(S^1)},\quad \forall \psi_1,\psi_2\in W_1(0).
\eeq
Denote

\[
\xi=(\xi_1,\xi_2,\cdots,\xi_m),\quad \Pi \psi=\sum_{j=1}^m\xi_j\varphi_j,
\]
such that $\sum_{j=1}^m\xi_j\varphi_j \in W_2(0)$ when $|\xi|$ is sufficiently small.
We now define $\Gamma:\,\R^m\mapsto \R$ as follows:

\beq\label{71}
\Gamma(\xi)=E\big(\mathscr{T}(\sum_{j=1}^m\xi_j\varphi_j)+w\big).
\eeq
It is clear that $\Gamma(\xi)$ is analytic in a small neighborhood of the origin in $\R^m$.

Straightforward calculations show that

\beq\label{53}
\frac{\partial \Gamma}{\partial\xi_j}=DE\big(\mathscr{T}(\sum_{j=1}^m\xi_j\varphi_j)+w\big)\cdot D\mathscr{T}\big(\sum_{j=1}^m\xi_j\varphi_j\big) \varphi_j.
\eeq
We infer from \eqref{52} that

\[
DE(u)\cdot v=\int_{S^1}\big(u_{\theta}v_{\theta}-A^2u v+\frac{\lambda}{u^2}v\big)d\theta,\,\ u, v\in H^2(S^1),
\]
which implies that

\beq\label{54}
DE(u)\cdot v=\int_{S^1}\big(-u_{\theta\theta}-A^2u +\frac{\lambda}{u^2}\big)v d\theta,\,\ u, v\in H^2(S^1).
\eeq
Because $w$ is an equilibrium, it follows from  \eqref{53} and \eqref{54}   that $\xi=0$ is a critical point of $\Gamma(\xi).$ By $\|\varphi_j\|_{L^2(S^1)}=1$, we infer from \eqref{55}, \eqref{53}, and \eqref{54} that

\[
\bigg|\frac{\partial \Gamma}{\partial\xi_j}\bigg|\leq \big\|\mathscr{M}\big(\mathscr{T}(\sum_{j=1}^m\xi_j\varphi_j)\big)\big\|_{L^2(S^1)} \big\|D\mathscr{T}\big(\sum_{j=1}^m\xi_j\varphi_j\big)\big\|_{\mathcal {L}(L^2,H^2)}\big\|\varphi_j\big\|_{L^2}\leq C \big\|\mathscr{M}\big(\mathscr{T}(\sum_{j=1}^m\xi_j\varphi_j)\big)\big\|_{L^2(S^1)}.
\]
Because $\Pi \psi=\sum_{j=1}^m\xi_j\varphi_j$, we derive from the above that

\begin{equation}\label{58}
\big|\nabla \Gamma(\xi)\big|\leq  C \|\mathscr{M}(\mathscr{T}(\Pi \psi))\|_{L^2(S^1)}\\
\leq C\big(\|\mathscr{M}(\psi)\|_{L^2(S^1)}+\|\mathscr{M}(\mathscr{T}(\Pi \psi))-\mathscr{M}(\psi)\|_{L^2(S^1)}\big).
\end{equation}
Recall from $\mathscr{N}(\psi)=\mathscr{M}(\psi)+\Pi \psi$ and \eqref{56} that

\begin{equation}\label{59}
\begin{split}
\|\mathscr{M}(\psi_1)-\mathscr{M}(\psi_2)\|_{L^2(S^1)}\leq &\|\Pi \psi_1-\Pi \psi_2\|_{L^2(S^1)}+\|\mathscr{N}(\psi_1)-\mathscr{N}(\psi_2)\|_{L^2(S^1)}\\
\leq & C \|\psi_1- \psi_2\|_{H^2(S^1)}.
\end{split}
\end{equation}
Note also from \eqref{57} that

\beq \label{60}
\|\mathscr{T}(\Pi \psi)-\psi\|_{H^2(S^1)}=\|\mathscr{T}(\Pi \psi)-\mathscr{T}(\mathscr{N}(\psi))\|_{H^2(S^1)}
\leq C\|\Pi \psi-\mathscr{N}(\psi)\|_{L^2(S^1)}=C\|\mathscr{M}(\psi)\|_{L^2(S^1)}.
\eeq
We then deduce from \eqref{58}, \eqref{59}, and \eqref{60} that

\begin{equation}\label{511}
|\nabla \Gamma(\xi)|
\leq   C\|\mathscr{M}(\psi)\|_{L^2(S^1)}.
\end{equation}
If $v\in W_1(0),$ then $v+t(\mathscr{T}(\Pi v)-v)\in W_1(0)$ for $t\in[0,1]$. We thus  infer from \eqref{71}, \eqref{59}, and \eqref{60} that for $v=w+\psi,$

\begin{eqnarray}\label{512}
\nonumber|E(v)-\Gamma(\xi)|&=&|E(v)-E(\mathscr{T}(\Pi \psi)+w)|\\
\nonumber&=&\Big|\int_0^1\frac{d}{dt}E\big(v+(1-t)(\mathscr{T}(\Pi \psi)-\psi)\big)dt\Big|\\
\nonumber&=&\Big|\int_0^1DE\big(v+(1-t)(\mathscr{T}(\Pi \psi)-\psi)\big)\cdot \big(\mathscr{T}(\Pi \psi)-\psi\big)dt\Big|\\
\nonumber&\leq &\max_{0\leq t\leq 1}\|\mathscr{M}(\psi+(1-t)(\mathscr{T}(\Pi \psi)-\psi))\|_{L^2(S^1)}\|\mathscr{T}(\Pi \psi)-\psi\|_{L^2(S^1)}\\
\nonumber&\leq &\bigg(\|\mathscr{M}(\psi)\|_{L^2(S^1)}+C\|\mathscr{T}(\Pi \psi)-\psi\|_{H^2(S^1)}\bigg)\|\mathscr{T}(\Pi \psi)-\psi\|_{H^2(S^1)}\\
&\leq & C\|\mathscr{M}(\psi)\|_{L^2(S^1)}^2.
\end{eqnarray}
By the Lojasiewicz inequality (cf. \cite{L1}), there exist a small constant $\sigma>0$ and $\bbtheta\in(0,{1\over2})$ such that

\beq \label{510}
|\nabla \Gamma(\xi)|\geq |\Gamma(\xi)-\Gamma(0)|^{1-\bbtheta}\quad \text{for} \,\ |\xi|\leq\sigma.
\eeq
Note that we can choose $\sigma>0$ small enough that   $\Pi \psi=\sum_{j=1}^m\xi_j\varphi_j\in W_2(0)$ holds for $|\xi|\leq \sigma$. By the definition of $\Gamma(\xi)$, we   infer from \eqref{510} that if $\sigma>0$ is small enough, then

\beq \label{510T}
|\nabla \Gamma(\xi)|\geq |\Gamma(\xi)-E(w)|^{1-\bbtheta}\quad \text{for} \,\ |\xi|\leq\sigma.
\eeq
Applying the elementary inequality$|a+b|^{1-\bbtheta}\geq {1\over2} |a|^{1-\bbtheta}-{1\over2} |b|^{1-\bbtheta},$
we deduce from \eqref{511}, \eqref{512}, and \eqref{510T} that for $\psi\in W_1(0)$,

\beq\label{72}
C\|\mathscr{M}(\psi)\|_{L^2(S^1)}\geq |\nabla \Gamma(\xi)|
\geq  |\Gamma(\xi)-E(w)|^{1-\bbtheta}
\geq  {1\over2} |E(v)-E(w)|^{1-\bbtheta}-{1\over2}C^{1-\bbtheta} \|\mathscr{M}(\psi)\|_{L^2(S^1)}^{2(1-\bbtheta)}.
\eeq
Because $0<\bbtheta<{1\over2},$ we have $2(1-\bbtheta)>1.$ It then follows from \eqref{72} that

\beq \label{515}
\|\mathscr{M}(\psi)\|_{L^2(S^1)}\geq C |E(v)-E(w)|^{1-\bbtheta},\, \quad \psi\in W_1(0).
\eeq
Choose both $\varepsilon>0$ and  $\sigma>0$   sufficiently small, such that

 \beq\label{514}
 C|E(v)-E(w)|^{-\varepsilon}\geq 1,\,\ \mbox{if}\,\ \|\psi\|_{H^2(S^1)}\leq\sigma.
 \eeq
Combining \eqref{515} with \eqref{514}, there exist sufficiently small constants $\varepsilon>0$ and  $\sigma>0$ such that

\[
\|\mathscr{M}(\psi)\|_{L^2(S^1)}\geq  |E(v)-E(w)|^{1-\bbtheta'}
\]
holds for $\|\psi\|_{H^2(S^1)}\leq\sigma$, where $0<\bbtheta':=\bbtheta-\varepsilon<{1\over2}.$
This therefore completes the proof for the case where $\dim(ker L)=m>0$.

As for the case where $\dim(ker L)=0,$ similar to \eqref{512}, one can also obtain that there exists $\sigma>0$ depending on $w$ such that

\[
|E(v)-E(w)|\leq C \|\mathscr{M}(\psi)\|^2_{L^2(S^1)}, \quad \|v-w\|_{H^2(S^1)}=\|\psi\|_{H^2(S^1)}<\sigma,
\]
which can be rewritten as

\[
\|\mathscr{M}(\psi)\|_{L^2(S^1)}\geq C |E(v)-E(w)|^{1-\frac{1}{2}}.
\]
Then \eqref{519} follows by an analysis similar to that above.
This completes the proof of the lemma.
\qed
 \end{appendices}

\section*{Acknowledgement}
The authors would like to thanks the anonymous referees for many valuable comments and suggestions which lead to some improvements of the present paper.


\end{document}